\documentstyle[12pt]{article}
\newcommand{\mn}{\sl}
\def\afterthmseparator{}
\makeatletter
\renewcommand{\@begintheorem}[2]{\trivlist
      \item[\hskip \labelsep{\bf #1\ #2\unskip\afterthmseparator}]\mn}
\renewcommand{\@opargbegintheorem}[3]{\trivlist
      \item[\hskip \labelsep{\bf #1\ #2\ (#3)\unskip\afterthmseparator}]\mn}
\makeatother
\newtheorem{theorem}{Theorem}[section]
\newtheorem{lemma}[theorem]{Lemma}

\newtheorem{proposition}[theorem]{Proposition}
\newtheorem{rem}[theorem]{Remark}
\newenvironment{remark}{\renewcommand{\mn}{\rm} \begin{rem}}{\end{rem}}
\newtheorem{probl}[theorem]{Problem}

\newtheorem{df}[theorem]{Definition}
\newenvironment{definition}{\renewcommand{\mn}{\rm} \begin{df}}{\end{df}}
\newtheorem{exmpl}[theorem]{Example}

\newcounter{oq}
\newcommand{\que}{\refstepcounter{oq}\par{\bf \theoq.}~}

\newcommand{\qed}{$\;\;\;\Box$}
\newcommand{\bull}{$\;\;\;$\vrule height .9ex width .8ex depth -.1ex}
\newcounter{proof}
\newenvironment{proof}{%
\stepcounter{proof}%
\par\smallbreak\noindent{\bf Proof.~}}%
{\unskip\nobreak\hfill \bull \par\medbreak}
\newcommand{\noproof}{\unskip\nobreak\hfill \bull}
\newenvironment{subproof}{\par\noindent{\sl Proof of Claim.~}}%
{\qed}
\newcounter{claim}[proof]
\renewcommand{\theclaim}{\Alph{claim}}
\newenvironment{claim}{\refstepcounter{claim}%
\par\medskip\par\noindent{\it Claim~\theclaim.~}~\rm}%
{\par\smallskip\par}
\newcommand{\case}[2]{\par{\it Case #1:}\/ #2}
\newcommand{\subcase}[2]{\par{\it Subcase #1:}\/ #2}

\newcommand{\function}[2]{:#1 \rightarrow #2}
\newcommand{\of}[1]{\left( #1 \right)}
\newcommand{\setdef}[2]{\left\{\hspace{0.5mm}#1:\hspace{0.5mm} #2\right\}}
\newcommand{\setdeff}[2]{\left\{\hspace{0.5mm}#1%
\hspace{0.5mm} \Bigl| \hspace{0.5mm} #2\right\}}
\newcommand{\ant}{\hspace{-0.6em}}    

\newcommand{\gd}{\delta}
\newcommand{\gs}{\sigma}

\newcommand{\And}{\wedge}
\newcommand{\Or}{\vee}

\newcommand{\hide}[1]{}
\newcommand{\refeq}[1]{(\ref{#1})}

\newcommand{\EF}{Ehrenfeucht}
\newcommand{\BS}{Bernays-Sch\"onfinkel}
\newcommand{\WL}{Weisfeiler-Lehman}
\newcommand{\game}{\mbox{\sc Ehr}}
\newcommand{\D}[1]{\mbox{\rm D\hspace{2pt}}(#1)}
\newcommand{\DD}[2]{\mbox{\rm D}_{#1} (#2)}
\newcommand{\compl}[1]{\overline{#1}}
\newcommand{\Nest}{\mathop{\rm Nest}\nolimits}
\newcommand{\nest}[1]{\Nest(#1)}
\newcommand{\Qr}{\mathop{\rm qr}\nolimits}
\newcommand{\qr}[1]{\Qr(#1)}
\newcommand{\Alt}{\mathop{\rm alt}\nolimits}
\newcommand{\alt}[1]{\Alt(#1)}
\newcommand{\phicong}{\mathbin{{\cong}_\phi}}
\newcommand{\xeq}{\mathbin{{\equiv}_X}}
\newcommand{\xcupeq}{\mathbin{{\equiv}_{X\cup\{u\}}}}
\newcommand{\xxeq}{\mathbin{{\equiv}_{X'}}}
\newcommand{\xkeq}{\mathbin{{\equiv}_{X_k}}}
\newcommand{\xkkeq}{\mathbin{{\equiv}_{X_{k+1}}}}
\newcommand{\xjeq}{\mathbin{{\equiv}_{X_j}}}
\newcommand{\xjjeq}{\mathbin{{\equiv}_{X'_j}}}
\newcommand{\notxjeq}{\mathbin{{\not\equiv}_{X'_j}}}
\newcommand{\phieq}{\mathbin{{\equiv}_\phi}}
\newcommand{\phimeq}{\mathbin{{\equiv}_{\phi_m}}}
\newcommand{\phijeq}{\mathbin{{\equiv}_{\phi_j}}}
\newcommand{\phikeq}{\mathbin{{\equiv}_{\phi_k}}}
\newcommand{\phiimineq}{\mathbin{{\equiv}_{\phi_{i-1}}}}
\newcommand{\notphiieq}{\mathbin{{\not\equiv}_{\phi_i}}}
\newcommand{\notxjjeq}{\mathbin{{\not\equiv}_{X'_j}}}
\newcommand{\notphijeq}{\mathbin{{\not\equiv}_{\phi_j}}}
\newcommand{\notphimeq}{\mathbin{{\not\equiv}_{\phi_m}}}
\newcommand{\xapp}{\mathbin{{\approx}_X}}
\newcommand{\notxjapp}{\mathbin{{\not\approx}_{X'_j}}}
\newcommand{\id}{\mbox{\rm id}}
\newcommand{\calC}{{\cal C}}

\newcommand{\dom}{\mathop{\rm dom}\nolimits}
\newcommand{\range}{\mathop{\rm range}\nolimits}
\newcommand{\ofcurle}[1]{\left\{ #1 \right\}}
\newcommand{\PR}[2]{\mbox{\rm P}_{#1} (#2)}
\newcommand{\bs}[1]{\mbox{\rm BS\hspace{2pt}}(#1)}
\newcommand{\Bs}[2]{\mbox{\rm BS}_{#1} (#2)}
\newcommand{\Dist}{\mbox{\rm Dist}}
\newcommand{\Iso}{\mbox{\rm Iso}}
\newcommand{\notaeq}{\mathbin{{\not\equiv}_{\compl A}}}
\newcommand{\beq}{\mathbin{{\equiv}_{B}}}

\newcommand{\dfn}[1]{\D{#1}}
\newcommand{\Dfn}[2]{\DD{#1}{#2}}
\newcommand{\I}[1]{\mbox{\rm I\hspace{2pt}}(#1)}
\newcommand{\II}[2]{\mbox{\rm I}_{#1} (#2)}
\newcommand{\add}{\oplus}
\newcommand{\baru}{{\bar u}}
\newcommand{\tiu}{{\tilde u}}
\newcommand{\hatu}{{\hat u}}
\newcommand{\hatU}{{\hat U}}

\newif\ifnotesw\noteswtrue

\title{
Descriptive Complexity of Finite Structures:\\
Saving the Quantifier Rank
}

\author{Oleg Pikhurko%
\thanks{Department of Mathematical Sciences,
Carnegie Mellon University, Pittsburgh, PA 15213-3890.
Web: {\tt http://www.math.cmu.edu/\~{}pikhurko/}}
\ and Oleg Verbitsky%
\thanks{Dept.~of Mechanics and Mathematics, Kyiv University, Ukraine.
E-mail: {\tt oleg@ov.litech.net}
}}

\date{}

\begin{document}
\maketitle

\begin{abstract}
We say that a first order formula $\Phi$ {\em distinguishes\/} a structure
$M$ over a vocabulary $L$ from another structure $M'$
over the same vocabulary
if $\Phi$ is true on $M$ but false on $M'$. A formula $\Phi$
{\em defines\/} an $L$-structure $M$ if $\Phi$ distinguishes $M$ from
any other non-isomorphic $L$-structure $M'$. A formula $\Phi$
{\em identifies\/} an $n$-element $L$-structure $M$ if $\Phi$
distinguishes $M$ from any other non-isomorphic $n$-element
$L$-structure~$M'$.

We prove that every $n$-element structure $M$ is identifiable by
a formula with quantifier rank less than $(1-\frac1{2k})n+k^2-k+2$
and at most one quantifier alternation, where $k$ is the maximum
relation arity of $M$. Moreover, if the automorphism group of $M$
contains no transposition of two elements, the same result holds
for definability rather than identification.

The \BS\/ class consists of prenex formulas in which the existential
quantifiers all precede the universal quantifiers.
We prove that every $n$-element structure $M$ is identifiable
by a formula in the \BS\/ class with less than $(1-\frac1{2k^2+2})n+k$
quantifiers. If in this class of identifying formulas we restrict the
number of universal quantifiers to $k$, then less than $n-\sqrt n+k^2+k$
quantifiers suffice to identify $M$ and, as long as we keep the
number of universal quantifiers bounded by a constant, at total
$n-O(\sqrt n)$ quantifiers are necessary.
\end{abstract}

\clearpage

\tableofcontents

\clearpage

\section{Introduction}

Let $M$ be a structure over a vocabulary $L$. A closed first order formula
$\Phi$ with relation symbols in $L\cup\{=\}$ is either true or false
on $M$. If $M'$ is another $L$-structure isomorphic with $M$, then
$\Phi$ is equally true or false on $M$ and $M'$. On the other hand,
if $M$ is finite and $M'$ is non-isomorphic with $M$, then there is
a formula $\Phi_{M,M'}$ that is true on $M$ and false on $M'$.
As it is well known, for infinite structures this is not necessary true.
In this paper, however, we deal only with finite structures.
We call the number of elements of a structure $M$ its {\em order}.

If a first order formula $\Phi$ is true on $M$ but false on $M'$,
we say that $\Phi$ {\em distinguishes $M$ from $M'$}. We say that
$\Phi$ {\em defines\/} an $L$-structure $M$ if $\Phi$ distinguishes $M$
from any other non-isomorphic $L$-structure $M'$. Furthermore, a formula
$\Phi$ {\em identifies\/} a finite $L$-structure $M$ if $\Phi$
distinguishes $M$ from any other non-isomorphic $L$-structure $M'$
of the same order.

We address the question how simple a formula identifying
(defining) a finite structure can be. The complexity measure of a first
order formula we use here is the quantifier rank, that is, the maximum
number of nested quantifiers in a formula. Let $\I M$ (resp.\ $\dfn M$)
denote the minimum quantifier rank of a formula identifying (resp.\
defining) a structure $M$. We will pay a special attention to
formulas of restricted logical structure. The {\em alternation number\/}
of a formula $\Phi$ is the maximum number of quantifier alternations over all
possible sequences of nested quantifiers under the assumption that
$\Phi$ is reduced to its negation normal form, i.e., all negations
are assumed to occur only in front of atomic subformulas.
By $\II lM$ and $\Dfn lM$
we denote the variants of $\I M$ and $\dfn M$ for the class of formulas
with alternation number at most $l$.

We will estimate $\I M$ and $\dfn M$ as functions of the order of $M$.
The latter is denoted throughout the paper by $n$. A simple upper
bound for $\I M$ is
$$
\II0M\le n.
$$
Indeed, every structure $M$ is identified by formula
\begin{equation}\label{eq:id}
\exists x_1\ldots\exists x_n\of{
\bigwedge_{1\le i<j\le n}x_i\ne x_j\And \Psi_M(x_1,\ldots,x_n)
},
\end{equation}
where $\Psi_M$ is the conjunction that gives an account of all
relations between elements of $M$ and negations thereof.
For example, if $M$ consists of a single binary relation $R^M$ on the set
$\{1,\ldots,n\}$, then
$$
\Psi_M=
\bigwedge_{(i,j)\in R^M}R(x_i,x_j)\And
\bigwedge_{(i,j)\notin R^M}\neg R(x_i,x_j).
$$

It is an easy exercise to show that, if $M$ has only unary relations,
then $\II0M\le(n+1)/2$. In \cite{PVV} we prove the following results.
If $M$ has only unary and binary relations, then $\II1M\le(n+3)/2$.
In the particular case that $M$ is an ordinary undirected graph,
we are able to improve on the alternation number
by showing that then $\II0M\le(n+5)/2$.
It is not hard to show that these bounds are tight up to
a small additive constant.
If $M$ is a $k$-uniform hypergraph, we have the bound
$\II1M\le(1-1/k)n+2k-1$.

Here we continue the research initiated in \cite{PVV} and prove
a general upper bound
\begin{equation}\label{eq:id_general}
\II1M<\of{1-\frac1{2k}}n+k^2-k+2,
\end{equation}
where $k$, here and throughout, denotes the maximum relation arity of
the vocabulary~$L$.

A simple upper bound for $\dfn M$ is
$$
\Dfn1M\le n+1.
$$
An appropriate defining formula is
$$
\exists x_1\ldots\exists x_n\forall x_{n+1}\of{
\bigwedge_{1\le i<j\le n}(x_i\ne x_j)\ \ \And\ \
\bigvee_{i=1}^n (x_{n+1}=x_i)\ \
\And\ \ \Psi_M(x_1,\ldots,x_n)
},
$$
where $\Psi_M$ is as in \refeq{eq:id}. The upper bound of $n+1$
is generally best possible. For example, we have $\dfn{M_n}=n+1$
if $M_n$ consists of the single totally true unary relation or
is a complete graph on $n$ vertices. However, for a quite representative
class of structures we are able to prove a better bound. We call
a structure {\em irredundant\/} if its automorphism group contains
no transposition of two elements. Similarly to \refeq{eq:id_general},
for any irredundant structure $M$ we obtain
\begin{equation}\label{eq:def_irred}
\Dfn1M<\of{1-\frac1{2k}}n+k^2-k+1.
\end{equation}
This is a qualitative extension of a result in \cite{PVV}, where the
bound $\Dfn1M\le n/2+2$ is proved for any irredundant structure $M$
with maximum relation arity 2.
On the other hand, there are simple examples of irredundant structures
with $\dfn M\ge n/4$ (see Remark \ref{rem:irred}).

In fact, the bound $\Dfn1M<(1-\frac1{2k})n+k^2-k+2$
may not hold only for structures with a simple, easily recognizable
property. Namely, given elements $u$ and $v$ of $M$, let us call
them {\em similar\/} if the transposition of $u$ and $v$ is an
automorphism of $M$. It turns out that, either we have the upper
bound for $\Dfn1M$ or otherwise $M$ has more than
$(1-\frac1{2k})n+(k-1)^2$ pairwise similar elements.
In the latter case we are able to easily compute the value
of $\dfn M$ up to an additive constant of $k$.
For graphs such a dichotomy result was obtained in~\cite{PVV}.

Furthermore, we address the identification of finite structures
by formulas of the
simplest logical structure, namely, those in the {\em prenex normal form}
(or {\em prenex formulas}).
In this case the quantifier rank is just the number of quantifiers
occurring in a formula.
Let $\Sigma_1$ (resp.\ $\Pi_1$) consist of the existential
(resp.\ universal) prenex formulas. Furthermore,
let $\Sigma_i$ (resp.\ $\Pi_i$) be the extension of
$\Sigma_{i-1}\cup\Pi_{i-1}$ with prenex formulas whose
quantifier prefix begins with $\exists$ (resp.\ with $\forall$)
and has less than $i$ quantifier alternations. In particular,
$\Sigma_2$ is the well-known \BS\/ class of formulas
(see \cite{Fag} for the role of this class in finite model theory).
Define $\PR iM$ to be the minimum number of quantifiers in
a $\Sigma_i\cup\Pi_i$ formula identifying a structure $M$.
Similarly, let $\bs M$ be the minimum number of quantifiers of
an identifying formula in the \BS\/ class $\Sigma_2$.
We hence have the following hierarchy:
\begin{equation}\label{eq:hierarch}
\begin{array}{cl}
\I M\le\II{i-1}M\le\PR{i}M\le\PR{i-1}M,&\quad i\ge1;\\[2mm]
\PR{2}M\le\bs M\le\PR{1}M\le n.&
\end{array}
\end{equation}
The upper bound of $n$ is here due to the identifying formula
\refeq{eq:id_general}. The bound $\PR1M\le n$ is generally best possible.
It is attained, for example, if $M$ consists of the single unary
relation true on all but one elements of the structure.

Our concern becomes therefore $\bs M$, the next member at the top
of the hierarchy \refeq{eq:hierarch}. We prove that
\begin{equation}\label{eq:bs}
\bs M < \of{1-\frac1{2k^2+2}}n+k.
\end{equation}
Though the multiplicative constant in \refeq{eq:bs} is worse than
that in the bound \refeq{eq:id_general}, the bound \refeq{eq:bs}
may be regarded as a qualitative strengthening of \refeq{eq:id_general}
because the class of formulas in the former result is much more limited
than that in the latter result. Curiously, the bound \refeq{eq:bs}
strengthens the bound \refeq{eq:id_general} also quantitatively if we
consider a somewhat unusual complexity measure of a formula, namely,
the total number of quantifiers occurring in it.

If we restrict the number of universal quantifiers to a constant,
Bernays-Sch\"onfin\-kel formulas become much less powerful. Let
$\Bs qM$ denote the minimum total number of quantifiers in a
\BS\/ formula identifying $M$ with at most $q$ universal quantifiers.
We prove that $\Bs kM<n-\sqrt n+k^2+k$ and that
$\Bs qM\ge n-O(\sqrt n)$ as long as $q$ is bounded by a constant.

To prove \refeq{eq:id_general}, we use the characterization of
the quantifier rank of a formula distinguishing structures $M$ and
$M'$ as the length of the {\em \EF\/ game\/} on $M$ and $M'$ \cite{Ehr}
(an essentially equivalent characterization in terms of partial isomorphisms
between $M$ and $M'$ and extensions thereof is due to Fra\"\i ss\'e
\cite{Fra}). Unlike \refeq{eq:id_general}, our proof of \refeq{eq:bs}
uses a direct approach. Nevertheless, both the results share the
same background which is based on the notion of a base of a structure~$M$.

Given a set $X$ of elements of $M$ and elements $u$ and $v$ of $M$,
we say that $X$ {\em separates\/} $u$ and $v$ if the extension of
the identity map of $X$ onto itself taking $u$ to $v$ is not a
partial automorphism of $M$.
Clearly, no $X$ can separate similar $u$ and $v$. On the other hand,
if $X$ separates every two non-similar elements in the complement
of $X$, we call $X$ a {\em base\/} of $M$. Every $M$ trivially
has $(n-1)$-element bases. Our technical results imply that
a considerably smaller base always exists.\footnote{%
In fact, we do not state this explicitly. However, it is easy
to derive from the estimate \refeq{eq:X=n_le} that every structure
has a base with less than $(1-\frac1{2k^2+1})n$ elements.
On the other hand, there are structures whose all bases have at least
$n/2$ elements. A simple example is given by the graph with $m$
pairwise non-adjacent edges.}

\paragraph{Related work.}
Our paper is focused on the descriptive complexity of {\em individual\/}
structures as opposed to the descriptive complexity of {\em classes\/}
of structures. The latter is the subject of a large research area,
which is emphasized much on the monadic second order logic
(we refer the reader to the survey \cite{Fag} and textbooks \cite{EFl,Imm}).

The identification of graphs in first order logic is studied in
\cite{IKo,ILa,CFI,Gro1,Gro2} in aspects relevant to computer science.
The main focus of this line of research is on the minimum number
of variables in an identifying formula, where formulas are in
the first order language enriched by {\em counting\/} quantifiers.
This complexity measure of a formula corresponds to the dimension
of the \WL\/ algorithm that succeeds in finding a canonic form
of a graph~\cite{CFI}.

The present paper studies, in a sense, the worst case descriptive
complexity of a structure. Two other possibilities,
the ``best'' and average structures, are considered
in \cite{PSV} and \cite{KPSV} in the case of graphs.

\paragraph{Organization of the paper.}
In Section \ref{s:back} we explain the notation used throughout
the paper, recall some basic definitions, define the \EF\/ game
and state its connection to distinguishing
non-isomorphic structures in first order logic. In Section \ref{s:base}
we introduce some relations, partitions,
transformations, and constructions over a finite structure and explore
their properties. The main task performed in this section is construction
of a particular base in an arbitrary structure. We will benefit from
these preliminaries while proving our both main results, bounds
\refeq{eq:id_general} and \refeq{eq:bs}, in Sections \ref{s:id} and
\ref{s:bs} respectively.
In Section \ref{s:id} we also prove the bound \refeq{eq:def_irred}
and the other definability results.
Section \ref{s:unibound} is devoted to identification
by \BS\/ formulas with bounded number of universal quantifiers.
In Section \ref{s:graphs} we focus on graphs and improve the
bound \refeq{eq:bs} for this class of structures. We conclude with
a list of open problems in Section~\ref{s:open}.

\section{Background}\label{s:back}

\subsection{Notation}\label{s:notation}

Writing $\bar u\in U^k$ for a set $U$ and a positive integer $k$,
we mean that $\bar u=(u_1,\ldots,u_k)$ with $u_i\in U$ for every
$i\le k$. If $u,v\in U$, then $\bar u^{(uv)}$ denotes\footnote{%
The double use of the character $u$ here should not be confusing:
We will often use $u$ to denote a single element of a sequence~$\bar u$.
}
the result
of substituting $v$ in place of every occurrence of $u$ in $\bar u$
and substituting $u$ in place of every occurrence of $v$ in $\bar u$.
Here $(uv)$ denotes the transposition of $u$ and $v$, that is,
the permutation of $U$ interchanging $u$ and $v$ and leaving
the remaining elements unchanged.
Given a function $\phi$ defined on $U$, we extend it over $U^k$
by $\phi(\bar u)=(\phi(u_1),\ldots,\phi(u_k))$ for $\bar u\in U^k$.

Notation $\id_U$ stands for the identity map of a set $U$ onto itself.
The domain and range of a function $f$ are denoted by $\dom f$
and $\range f$ respectively.
By $f^{(k)}$ we denote the $k$-fold composition of~$f$.

\subsection{Basic definitions}\label{s:defn}

A {\em $k$-ary relation $R$ on a set $V$\/} (or a {\em relation $R$
of arity $k$\/}) is a function from $V^k$ to $\{0,1\}$.
A {\em vocabulary\/} is a finite sequence $R_1,\ldots,R_m$ of relation
symbols along with a sequence $k_1,\ldots,k_m$ of positive integers,
where each $k_i$ is the arity of the respective $R_i$.
If $L$ is a vocabulary, a finite {\em structure $A$ over $L$\/}
(or an {\em $L$-structure $A$\/}) is a finite set $V(A)$, called the
{\em universe\/}, along with relations $R_1^A,\ldots,R_m^A$, where
$R_i^A$ has arity $k_i$. The {\em order\/} of $A$ is the number
of elements in the universe $V(A)$. If $U\subseteq V(A)$, then $A$
induces on $U$ the structure $A[U]$ with the universe $V(A[U])=U$
and relations $R^{A[U]}_1,\ldots,R^{A[U]}_m$ such that
$R_i^{A[U]}\bar a=R_i^A\bar a$ for every $\bar a\in U^{k_i}$.
Two $L$-structures $A$ and $B$ are {\em isomorphic\/} if there is
a one-to-one map $\phi\function{V(A)}{V(B)}$, called {\em an isomorphism
from $A$ to $B$\/}, such that
$R_i^A\bar a=R_i^B\phi\bar a$ for every $i\le m$ and all
$\bar a\in V(A)^{k_i}$.
An {\em automorphism of $A$\/} is an isomorphism from $A$ to itself.
If $U\subseteq V(A)$ and $W\subseteq V(B)$,
we call a one-to-one map $\phi\function UW$ a {\em partial isomorphism
from $A$ to $B$\/} if it is an isomorphism from $A[U]$ to $B[W]$.

Without loss of generality we assume first order formulas to be over
the set of connectives $\{\neg,\And,\Or\}$.

\begin{definition}
A {\em sequence of quantifiers\/} is a finite word over the alphabet
$\{\exists,\forall\}$. If $S$ is a set of such sequences, then
$\exists S$ (resp.\ $\forall S$) means the set of concatenations
$\exists s$ (resp.\ $\forall s$) for all $s\in S$. If $s$ is a sequence
of quantifiers, then $\bar s$ denotes the result of replacement of all
occurrences of $\exists$ to $\forall$ and vice versa in $s$. The
 set $\bar S$
consists of all $\bar s$ for $s\in S$.

Given a first order formula $\Phi$, its set of {\em sequences of nested
quantifiers\/} is denoted by $\nest\Phi$ and defined by induction as
follows:
\begin{enumerate}
\item
$\nest\Phi=\{\lambda\}$ if $\Phi$ is atomic, where $\lambda$ denotes the
empty word.
\item
$\nest{\neg\Phi}=\overline{\nest\Phi}$.
\item
$\nest{\Phi\And\Psi}=\nest{\Phi\Or\Psi}=\nest\Phi\cup\nest\Psi$.
\item
$\nest{\exists x\Phi}=\exists\nest\Phi$ and
$\nest{\forall x\Phi}=\forall\nest\Phi$.
\end{enumerate}

The {\em quantifier rank\/} of a formula $\Phi$, denoted by $\qr\Phi$,
is the maximum length of a string in $\nest\Phi$.

Given a sequence of quantifiers $s$, let $\alt s$ denote the number
of occurrences of $\exists\forall$ and $\forall \exists$ in $s$.
The {\em alternation number\/} of a first order formula $\Phi$
is the maximum $\alt s$ over $s\in\nest\Phi$.
\end{definition}

Given an $L$-structure $A$ and a closed first order formula $\Phi$
whose relation symbols are from $L\cup\{{=}\}$,
we write $A\models\Phi$ if $\Phi$ is true
on $A$ and $A\not\models\Phi$ otherwise.
Given $A$, a formula $\Psi(x_1,\ldots,x_m)$ with
$m$ free variables $x_1,\ldots,x_m$, and
a sequence $a_1,\ldots,a_m$ of elements in $V(A)$,
we write $A,a_1,\ldots,a_m\models\Psi(x_1,\ldots,x_m)$
if $\Psi(x_1,\ldots,x_m)$ is true on $A$ with each $x_i$
assigned the respective $a_i$.

If $B$ is another $L$-structure,
we say that a formula $\Phi$ {\em distinguishes $A$ from $B$\/} if
$A\models\Phi$ but $B\not\models\Phi$.
We say that $\Phi$ {\em defines\/} an $L$-structure $A$ (up to an isomorphism)
if $\Phi$ distinguishes $A$ from any non-isomorphic $L$-structure $B$.
We say that $\Phi$ {\em identifies\/} an $L$-structure $A$ of order $n$
(up to an isomorphism in the class of $L$-structures of the same order)
if $\Phi$ distinguishes $A$ from any non-isomorphic $L$-structure $B$
of order~$n$.

By $\D{A,B}$ (resp.\ $\DD l{A,B}$) we denote the minimum quantifier rank of
a formula (resp.\ with alternation number at most $l$) distinguishing
a structure $A$ from a structure $B$.
By $\dfn A$ (resp.\ $\Dfn lA$) we denote the minimum
quantifier rank of a formula defining $A$
(resp.\ with alternation number at most $l$).
By $\I A$ (resp.\ $\II lA$) we denote the minimum
quantifier rank of a formula identifying $A$
(resp.\ with alternation number at most $l$).

\begin{lemma}\label{lem:ddd}
Let $A$ be a finite structure over vocabulary $L$. Then the following
equalities hold true:
\begin{eqnarray*}
\dfn A&=&\max\setdef{\D{A,B}}{B\not\cong A},\hspace{60mm}\mbox{}\\
\Dfn lA&=&\max\setdef{\DD l{A,B}}{B\not\cong A},\\
\I A&=&\max\setdef{\D{A,B}}{B\not\cong A,\, |V(B)|=|V(A)|},\\
\II lA&=&\max\setdef{\DD l{A,B}}{B\not\cong A,\, |V(B)|=|V(A)|},
\end{eqnarray*}
where $\cong$ denotes the isomorphism relation between $L$-structures.
\end{lemma}

\begin{proof}
We prove the first equality; The proof of the others is similar.
Given an $L$-structure $B$ non-isomorphic with $A$,
let $\Phi_{B}$ be a formula
of minimum quantifier rank distinguishing $A$ from $B$, that is,
$\qr{\Phi_{B}}=\D{A,B}$. Let $R=\max_{B}\qr{\Phi_{B}}$.
We have $\dfn A\ge R$ because $\dfn A\ge\D{A,B}$ for every $B$.
To prove the reverse inequality $\dfn A\le R$, notice that $A$
is defined by the formula $\Phi=\bigwedge_{B}\Phi_{B}$
whose quantifier rank is $R$. The only problem is that $\Phi$
is an infinite conjunction (a $FO_{\infty\omega}$-formula).
However, as it is well known, over a fixed finite vocabulary
there are only finitely many inequivalent first order formulas of
bounded quantifier rank (see e.g.\ \cite{CFI,EFl,Imm}).
We therefore can reduce $\Phi$ to a finite conjunction.
\end{proof}

\subsection{The \EF\/ game}

Let $A$ and $B$ be structures over the same vocabulary
with disjoint universes.
The $r$-round \EF\/ game on $A$ and $B$,
denoted by $\game_r(A,B)$, is played by
two players, Spoiler and Duplicator, with $r$ pairwise distinct
pebbles $p_1,\ldots,p_r$, each given in duplicate. Spoiler starts the game.
A {\em round\/} consists of a move of Spoiler followed by a move of
Duplicator. In the $s$-th round Spoiler selects one of
the structures $A$ or $B$ and places $p_s$ on an element of this structure.
In response Duplicator should place the other copy of $p_s$ on an element
of the other structure. It is allowed to place more than one pebble on the
same element. We will use $a_s$ (resp.\ $b_s$) to
denote the element of $A$ (resp.\ $B$) occupied by $p_s$, irrespectively
of who of the players places the pebble on this element.
If after every of $r$ rounds it is true that
$$
a_i=a_j\mbox{ iff } b_i=b_j\mbox{ for all }i,j\le s,
$$
and the component-wise correspondence between $(a_1,\ldots,a_s)$ and
$(b_1,\ldots,b_s)$ is a partial isomorphism from $A$ to $B$, this is
a win for Duplicator;  Otherwise the winner is Spoiler.

The {\em $l$-alternation\/} \EF\/ game on $A$ and $B$ is a variant
of the game in which Spoiler is allowed to switch from one structure
to another at most $l$ times during the game, i.e., in at most $l$
rounds he can choose the structure other than that in the preceding round.

The following statement provides us with a robust technical tool.

\begin{lemma}\label{lem:loggames}
Let $A$ and $B$ be non-isomorphic structures over the same vocabulary.

\begin{enumerate}
\item
$\D{A,B}$ equals the minimum $r$ such that Spoiler has a winning
strategy in $\game_r(A,B)$.
\item
$\DD l{A,B}$ equals the minimum $r$ such that Spoiler has a winning
strategy in the $l$-alternation $\game_r(A,B)$.\noproof
\end{enumerate}
\end{lemma}

\noindent
We refer the reader to \cite[Theorem 1.2.8]{EFl}, \cite[Theorem 6.10]{Imm},
or \cite[Theorem 2.3.1]{Spe} for the proof of the first
claim and to \cite{Pez} for the second claim.

\section{Exploring structural properties of finite structures}\label{s:base}

\subsection{A few useful relations}\label{ss:userel}

Throughout this section we are given an arbitrary finite structure $M$ over
vocabulary $L$. We abbreviate $V=V(M)$.

\begin{definition}
For $a,b\in V$ we write $a\sim b$ if the transposition $(ab)$ is
an automorphism of $M$. In other words, $a\sim b$ if, for every $l$-ary
relation $R$ of $M$, we have $R\bar a=R\bar a^{(ab)}$ for all
$\bar a\in V^l$.
\end{definition}

\begin{lemma}\label{lem:simiseq}
$\sim$ is an equivalence relation on~$V$.
\end{lemma}

\begin{proof}
The relation is obviously reflexive and symmetric. The transitivity
follows from the facts that the composition of automorphisms is
an automorphism and that the transposition $(ac)$ is decomposed
into a composition of $(ab)$ and~$(bc)$.
\end{proof}

Given $X\subset V$, we will denote its complement by $\compl X=V\setminus X$.

\begin{definition}
Let $X\subset V$ and $a,b\in\compl X$. We write $a\xeq b$ if
$\id_X$ extends to an isomorphism from $M[X\cup\{a\}]$ to
$M[X\cup\{b\}]$. In other words, for every $l$-ary
relation $R$ of $M$, we have $R\bar a=R\bar a^{(ab)}$ for all
$\bar a\in (X\cup\{a\})^l$.

Furthermore, we write $a\xapp b$ if the transposition $(a,b)$
is an automorphism of $M[X\cup\{a,b\}]$. In other words, for every $l$-ary
relation $R$ of $M$, we have $R\bar a=R\bar a^{(ab)}$ for all
$\bar a\in (X\cup\{a,b\})^l$.
\end{definition}

Clearly, $a\xapp b$ implies $a\xeq b$. It is also clear that
$\xeq$ is an equivalence relation on $\compl X$.
In contrast to this, simple examples show that $a\xapp b$ is generally
not an equivalence relation.

\begin{definition}\label{def:calc}
$\calC(X)$ is the partition of $\compl X$ into $\xeq$-equivalence classes.
Furthermore, $\calC^m(X)=\setdef{C\in\calC(X)}{|C|\le m}$.
\end{definition}

The following lemma points some trivial but important properties of the
partition $\calC(X)$.

\begin{lemma}\label{lem:new}
\mbox{}

\begin{enumerate}
\item
If $X_1\subseteq X_2$, then
$\calC(X_2)$ is a refinement of $\calC(X_1)$ on $\compl{X_2}$.
\item
For any $X$, the $\sim$-equivalence classes restricted to $\compl X$
refine the partition $\calC(X)$.\noproof
\end{enumerate}
\end{lemma}

In the sequel $M'$ denotes another $L$-structure.

\begin{definition}\label{def:phieq}
Let $\phi\function X{X'}$ be a partial isomorphism from $M$ to $M'$.
Let $a\in\compl X$ and $a'\in\compl{X'}$. We write $a\phieq a'$ if
$\phi$ extends to an isomorphism from $M[X\cup\{a\}]$ to
$M'[X'\cup\{a'\}]$.
\end{definition}

\begin{lemma}\label{lem:phieq}
Let $\phi\function X{X'}$ be a partial isomorphism from $M$ to $M'$.
Then the following claims are true.

\begin{enumerate}
\item
Assume that
$a\xeq b$ and $a'\xxeq b'$. Then $a\phieq a'$ iff
$b\phieq b'$.
\item
Assume that $a\phieq a'$ and $b\phieq b'$. Then $a\xeq b$ iff
$a'\xxeq b'$.
\item
Let $\bar\phi$ be a partial isomorphism from $M$ to $M'$ which is an
extension of $\phi$. If $a\in\dom\bar\phi\setminus X$,
then $a\phieq\bar\phi(a)$.
\item
Let $\bar\phi$ be a partial isomorphism from $M$ to $M'$ which is an
extension of $\phi$. Let $a,b\in\dom\bar\phi\setminus X$. Then $a\xeq b$ iff
$\bar\phi(a)\xxeq\bar\phi(b)$.\noproof
\end{enumerate}
\end{lemma}

The proof is easy. Item 1 of the lemma makes the following definition
correct.

\begin{definition}\label{def:phieqcl}
Let $\phi\function X{X'}$ be a partial isomorphism from $M$ to $M'$.
Let $C\in\calC(X)$ and $C'\in\calC(X')$. We write $C\phieq C'$ if
$a\phieq a'$ for some (equivalently, for all) $a\in C$ and $a'\in C'$.
\end{definition}

\subsection{A couple of useful transformations}\label{ss:transform}

Let $M$ be a finite structure of order $n$ with the maximum relation
arity $k$. Let $X\subseteq V(M)$.
We define two transformations that, if applicable to $X$, extend it
to a larger set.

\begin{description}
\item[{\em Transformation $T$.}]
If there exists a set $S\subseteq\compl X$ with at most $k-1$ elements
such that $|\calC(X\cup S)|>|\calC(X)|$, take the lexicographically first
such $S$ and set $T(X)=X\cup S$. Otherwise $T$ is not applicable to~$X$.
\item[{\em Transformation $E$.}]
Apply $T$ iteratively as long as it is applicable. The result is denoted
by $E(X)$. In other words, $E(X)=T^{(n)}(X)$. If $T$ is not applicable
at all, set $E(X)=X$.
\end{description}

\begin{lemma}\label{lem:cd}
Assume that $T$ is not applicable to $X$. If
$C\in\calC(X)\setminus\calC^2(X)$, then $a\xapp b$ for every $a,b\in C$.
\end{lemma}

\begin{proof}
Let $C\in\calC(X)$ and $|C|\ge3$. Given $a$ and $b$ in $C$, we have to
show that $a\xapp b$. In other words, our task is,
given an $l$-ary relation $R$ of $M$ and $\bar a\in (X\cup\{a,b\})^l$,
to show that $R\bar a=R\bar a^{(ab)}$. If $\bar a$ contains no occurrence
of $a$ or no occurrence of $b$, this equality is true because $a\xeq b$.
It remains to consider the case that $\bar a$ contains occurrences of both
$a$ and~$b$.

\begin{claim}\label{cl:uvw}
Let $u$, $v$, and $w$ be pairwise distinct elements in $C$.
Let $R$ be an $l$-ary relation of $M$ and $\bar u\in(X\cup\{u,v\})^l$
with occurrences of both $u$ and $v$. Then $R\bar u=R\bar u^{(vw)}$.
\end{claim}

\begin{subproof}
If $R\bar u\ne R\bar u^{(vw)}$, then removal of $u$ from $C$ to $X$
splits $C$ into at least two $\xcupeq$-subclasses, containing
$v$ and $w$ respectively. This contradicts the assumption that $T$
is not applicable to~$X$.
\end{subproof}

Let $c$ be an arbitrary element in $C\setminus\{a,b\}$.
Applying Claim \ref{cl:uvw} repeatedly three times, we obtain
$$
R\bar a=R\bar a^{(bc)}=R(\bar a^{(bc)})^{(ab)}=
R((\bar a^{(bc)})^{(ab)})^{(ac)}=R\bar a^{(bc)(ab)(ac)}=R\bar a^{(ab)},
$$
as required.
\end{proof}

\begin{lemma}\label{lem:expan}
$|E(X)\setminus X|\le(k-1)|\calC(E(X))\setminus\calC(X)|$.\noproof
\end{lemma}

\subsection{The many-layered base of a finite structure}

\begin{definition}\label{def:base}
Suppose that a finite structure $M$ with maximum relation arity $k$ is given.
For $X\subset V(M)$, let $Y(X)=\bigcup_{C\in\calC^{k+1}(X)}C$.
We set
\begin{eqnarray*}
&&X_0=Y_0=\emptyset,\\
&&X_i=E(X_{i-1}\cup Y_{i-1})\mbox{\ for\ }1\le i\le k,\\
&&Y_i=Y(X_i)\mbox{\ for\ }1\le i\le k,\\
&&X_{k+1}=X_k\cup Y_k,\\
&&Z=V(M)\setminus X_{k+1}.
\end{eqnarray*}
We will call $X_{k+1}$ {\em the base of\/} $M$.
\end{definition}

An important role of the base of a finite structure is due to the
following fact (cf.\ the more general Definition~\ref{def:base2}).

\begin{lemma}\label{lem:xeqsim}
On $Z$ the relations $\xkeq$, $\xkkeq$, and $\sim$ coincide.
\end{lemma}

\begin{proof}
We start with relations $\xkeq$ and $\sim$.
Assume on the contrary that $a\xkeq b$ but $a\not\sim b$ for some
$a,b\in Z$. The latter means that, for some $l$-ary relation $R$ of $M$ and
$\bar a\in V^l$ with at least one occurrence of $a$,
\begin{equation}\label{eq:Rab}
R\bar a^{(ab)}\ne R\bar a.
\end{equation}
Denote $A=\{a_1,\ldots,a_l\}\setminus\{a,b\}$.
Since $|A|\le k-1$ and the $Y_i$'s are pairwise disjoint,
there is $j\le k$ such that
\begin{equation}\label{eq:empty}
A\cap Y_j=\emptyset.
\end{equation}
Remove all elements from $A\setminus X_j$ to $X_j$ and set $X'_j=X_j\cup A$.
Due to \refeq{eq:Rab}, this operation has the effect that
\begin{equation}\label{eq:nonapp}
a\notxjapp b.
\end{equation}
No class in $\calC(X_j)$ can disappear completely: The classes in
$\calC^{k+1}(X_j)$ can only split up because of \refeq{eq:empty},
the classes in $\calC(X_j)\setminus\calC^{k+1}(X_j)$ can lose
up to $k-1$ elements and/or split up.

Since $a\xkeq b$ and $a,b\in Z$, both $a$ and $b$ belong to the same
$\xkeq$-class $C^*$ containing at least $k+2$ elements. Let $C$ be the
$\xjeq$-class such that $C^*\subseteq C$. We now show that $C$ is split
up after modifying $X_j$ and therefore $|\calC(X'_j)|>|\calC(X_j)|$,
making a contradiction to the construction of~$X_j$.

Indeed, if $a\notxjeq b$, we have two subclasses containing respectively
$a$ and $b$. If $a\xjjeq b$, it follows by Lemma \ref{lem:cd}
from \refeq{eq:nonapp} that the class in $\calC(X'_j)$ containing $a$
and $b$ is exactly $\{a,b\}$. After removing at most $k-1$ elements,
in $C$ there remain at least 3 elements and therefore $C$ must have
at least one more $\xjjeq$-subclass besides $\{a,b\}$.

Thus, on $Z$ the relations $\xkeq$ and $\sim$ are identical.
By Item 1 of Lemma \ref{lem:new}, on $Z$ the relation $\xkkeq$
refines $\xkeq$. By Item 2 of the same lemma the converse is also
true. It follows that on $Z$ the relations $\xkkeq$ and $\xkeq$
also coincide.
\end{proof}

\begin{lemma}\label{lem:unpebbl}
Let $n$ be the order of $M$ and $k$ be the maximum relation arity of $M$.
We have
\begin{equation}\label{eq:b0}
\sum^k_{i=1}|\calC^{k+1}(X_i)|+\frac{|Z|}2>\frac n{2k}+\frac12-\frac1{2k}
\mbox{\ \ if\ \ }k\ge 2
\end{equation}
and
\begin{equation}\label{eq:a0}
2k\sum^{k-1}_{i=1}|\calC^{k+1}(X_i)|+
(k+1)|\calC^{k+1}(X_k)|+(k-1)|\calC(X_k)|+|Z|\ge n+k-1.
\end{equation}
\end{lemma}

\begin{proof}
By Lemma \ref{lem:expan} we have
\begin{eqnarray}
|X_1|&\le&(k-1)(|\calC(X_1)|-1),\label{eq:x1}\\
|X_i\setminus(X_{i-1}\cup Y_{i-1})|&\le&
(k-1)(|\calC(X_i)|
-|\calC(X_{i-1}\cup Y_{i-1})|)\label{eq:xi}
\end{eqnarray}
for $2\le i\le k$. Note that
$$
|\calC(X_i)|=|\calC^{k+1}(X_i)|+|\calC(X_i)\setminus\calC^{k+1}(X_i)|
$$
and
$$
|\calC(X_i)\setminus\calC^{k+1}(X_i)|\le|\calC(X_i\cup Y_i)|
$$
for $1\le i\le k$.
The latter inequality is true because, according to Item 1 of
Lemma \ref{lem:new}, the partition $\calC(X_i\cup Y_i)$ is
a refinement of $\calC(X_i)\setminus\calC^{k+1}(X_i)$.
Combining it with \refeq{eq:x1} and \refeq{eq:xi}, we obtain
\begin{eqnarray}
|X_1|&\le&(k-1)(|\calC^{k+1}(X_1)|+|\calC(X_1\cup Y_1)|-1)\label{eq:x12}\\
|X_i\setminus(X_{i-1}\cup Y_{i-1})|&\le&
(k-1)(|\calC^{k+1}(X_i)|+\calC(X_i\cup Y_i)-|\calC(X_{i-1}\cup Y_{i-1})|).
\label{eq:xi2}
\end{eqnarray}
Summing up \refeq{eq:x12} and \refeq{eq:xi2} over all $2\le i\le k$, we have
\begin{equation}\label{eq:forpr}
|X_1|+\sum^k_{i=2}|X_i\setminus(X_{i-1}\cup Y_{i-1})|\le
(k-1)\of{\sum^k_{i=1}|\calC^{k+1}(X_i)|+|\calC(X_k\cup Y_k)|-1}.
\end{equation}
According to Lemma \ref{lem:xeqsim},
\begin{equation}\label{eq:a1}
\calC(X_k\cup Y_k)=\calC(X_k)\setminus\calC^{k+1}(X_k)
\end{equation}
and, as a consequence,
\begin{equation}\label{eq:b1}
|\calC(X_k\cup Y_k)|\le|Z|/(k+2).
\end{equation}
{}From \refeq{eq:forpr} we conclude, using \refeq{eq:a1}, that
\begin{equation}\label{eq:a2}
|X_1|+\sum^k_{i=2}|X_i\setminus(X_{i-1}\cup Y_{i-1})|\le
(k-1)\of{\sum_{i=1}^{k-1}|\calC^{k+1}(X_i)|+|\calC(X_k)|-1}
\end{equation}
and, using \refeq{eq:b1}, that
\begin{equation}\label{eq:b2}
|X_1|+\sum^k_{i=2}|X_i\setminus(X_{i-1}\cup Y_{i-1})|\le
(k-1)\of{\sum_{i=1}^{k}|\calC^{k+1}(X_i)|+\frac{|Z|}{k+2}-1}.
\end{equation}
Notice also a trivial inequality
\begin{equation}\label{eq:y}
|Y_i|\le(k+1)|\calC^{k+1}(X_i)|.
\end{equation}
It is easy to see that
\begin{equation}\label{eq:n}
n=|X_1|+\sum^k_{i=2}|X_i\setminus(X_{i-1}\cup Y_{i-1})|+
\sum^k_{i=1}|Y_i|+|Z|.
\end{equation}
Using \refeq{eq:a2} and \refeq{eq:y}, we derive from \refeq{eq:n} that
$$
n\le 2k\sum^{k-1}_{i=1}|\calC^{k+1}(X_i)|+
(k+1)|\calC^{k+1}(X_k)|+(k-1)|\calC(X_k)|+|Z|-(k-1),
$$
which implies \refeq{eq:a0}.
Using \refeq{eq:b2} and \refeq{eq:y}, we derive from \refeq{eq:n} that
\begin{equation}\label{eq:lastbut}
n\le 2k\sum^{k}_{i=1}|\calC^{k+1}(X_i)|+
\of{2-\frac3{k+2}}|Z|-(k-1),
\end{equation}
which implies \refeq{eq:b0}.
\end{proof}

\section{Identifying finite structures with smaller quantifier rank}%
\label{s:id}

\begin{theorem}\label{thm:main}
Let $L$ be a vocabulary with maximum relation arity $k$.
For every $L$-structure $M$ of order $n$ we have
$$
\II1{M} < \of{1-\frac1{2k}}n+k^2-k+2.
$$
\end{theorem}

The proof takes the next two subsections.
The case of $k=1$ is an easy exercise and we will assume that $k\ge2$.
According to Lemma \ref{lem:ddd}, it suffices to consider an arbitrary
$L$-structure $M'$ non-isomorphic with $M$ and of the same order $n$,
and estimate the value of $\DD1{M,M'}$.
We will design a strategy enabling Spoiler to win the \EF\/ game
on $M$ and $M'$ in less than $(1-\frac1{2k})n+k^2-k+2$ moves
with at most one alternation between the structures.
This will give us the desired bound by Lemma~\ref{lem:loggames}.

\subsection{Spoiler's strategy}\label{ss:strategy}

The strategy splits play into $k+2$ phases. Spoiler will play almost
all the time in $M$, possibly with one alternation from $M$ to $M'$
at the end of the game. For each vertex $v\in V(M)$ selected by Spoiler
up to Phase $i$, let $\phi^*_i(v)$ denote the vertex in $V(M')$ selected
in response by Duplicator. Thus, each subsequent $\phi^*_{i+1}$
extends $\phi^*_i$. Provided Phase $i$ has been already finished but
the game not yet, $\phi^*_i$ is a partial isomorphism from $M$ to $M'$.
Under the same condition, it will be always the case that
$\dom\phi^*_i\subseteq X_i$. We will use notation
$\tilde Y_{i-1}=\dom\phi^*_i\cap Y_{i-1}$.
Recall that the sets $X_i$ and $Y_i$ are defined by Definition \ref{def:base}
so that $Y_{i-1}\subset X_i$.

\medskip

{\sc Phase 1.}

Spoiler selects all vertices in $X_1$. Let $X'_1=\phi^*_1(X_1)$.

{\sc End of phase description.}

\medskip

{\sc Phase $j+1$, $1\le j\le k$.}

Our description of Phase $j+1$ is based on the assumption that Phase $j$
is complete but the game is not finished yet and that the following
conditions are true for every $1\le i\le j$.

\begin{description}
\item[{\em Condition 1.}]
$\phi^*_i$ has a unique extension $\phi_i$ over the whole $X_i$ that is
a partial isomorphism from $M$ to $M'$. Let $X'_i=\phi_i(X_i)$.
\item[{\em Condition 2.}]
There is a one-to-one correspondence between the partitions
$\calC^{k+1}(X_{i-1})$ and $\calC^{k+1}(X'_{i-1})$ such that,
if $C'\in\calC^{k+1}(X'_{i-1})$ corresponds to $C\in\calC^{k+1}(X_{i-1})$,
then $C\phiimineq C'$ and $|C|=|C'|$.
\item[{\em Condition 3.}]
For every $C\in\calC^{k+1}(X_{i-1})$, $\phi^*_i$ is defined on all
but one elements of $C$. Denote $\tilde C=\dom\phi^*_i\cap C$.
Then $\phi^*_i(\tilde C)\subset C'$, where $C'$ corresponds to $C$
according to Condition 2. Furthermore, $\phi_i$ takes the single
element in $C\setminus\tilde C$ to the single element in
$C'\setminus\phi^*_i(\tilde C)$. Thus, $\phi_i(C)=C'$.

For the further references we denote the set $\phi_i(Y_{i-1})=Y(X'_{i-1})$
by $Y'_{i-1}$ and its subset $\phi^*_i(\tilde Y_{i-1})$ by $\tilde Y'_{i-1}$.
\end{description}

Condition 1 is true for $i=1$ because $\phi^*_1$ is defined on the whole
$X_1$. For the sake of technical convenience, we set $X'_0=X_0=\emptyset$.
We suppose that $n>k+1$ (otherwise Theorem \ref{thm:main} is trivially true).
This implies that $\calC^{k+1}(X_0)=\calC^{k+1}(X'_0)=\emptyset$ and
makes Conditions 2 and 3 for $i=1$ trivially true. For $i>1$ Conditions 1--3
follow by induction from Claim \ref{cl:cond} below.

In the sequel we will intensively exploit the following notion.
We say that a pair $(a,a')\in V(M)\times V(M')$ is {\em $i$-threatening\/}
(for Duplicator) if $a$ and $a'$ are selected by the players in the
same round after Phase $i$ and
\begin{itemize}
\item
$a\notin X_i$ or $a'\notin X'_i$,
\item
$a\notphiieq a'$.
\end{itemize}
We now start description of the phase. It consists of two parts.

{\em Part 1.}
As long as no $i$-threatening pair arises for $1\le i\le j$, Spoiler
selects all but one elements in each class $C\in\calC^{k+1}(X_j)$.
The set of the vertices selected in $C$ will be denoted by $\tilde C$.
Furthermore,
Spoiler selects all vertices in
$X_{j+1}\setminus(X_j\cup Y_j)$. As soon as an $i$-threatening pair
for some $1\le i\le j$ arises, Spoiler switches to the strategy
given by Claim \ref{cl:threat} below and wins in at most $(i-1)(k-1)$ moves.

{\em Part 2.}
Assume that Part 1 finishes and Duplicator still does not lose.
Then, if Spoiler is able to win in at most $k$ next moves irrespective
of Duplicator's strategy, he does so and the game finishes.
If he is not able to win but able in at most $k$ moves
to enforce creating an $i$-threatening
pair for some $i\le j$, he does so and wins in at most $(i-1)(k-1)$
subsequent moves using the strategy of  Claim \ref{cl:threat}.
Otherwise Phase $j+1$ is complete and the next Phase $j+2$ starts.

{\sc End of phase description.}

\begin{claim}\label{cl:conserv}
Let $i\le k+1$. Suppose that Phase $i$ is finished and Conditions 1--3
are met for $i$ and all its preceding values. Assume that $a\in V(M)$ and
$a'\in V(M')$ are selected by the players in the same round after Phase $i$
and neither of them has been selected before. If
\begin{itemize}
\item
$a\in X_i$ but $a'\ne\phi_i(a)$

or
\item
$a'\in X'_i$ but $a\ne\phi_i^{-1}(a')$,
\end{itemize}
then the pair $(a,a')$ is $m$-threatening for some $m<i$.
\end{claim}

\begin{subproof}
Let $m$, $1\le m<i$, be the largest index such that neither $a\in X_m$
nor $a'\in X'_m$. Then $a\in X_{m+1}$ or $a'\in X'_{m+1}$. We consider
the former case (the analysis of the latter case is symmetric).
By Condition 3, $a\in Y_m\setminus\tilde Y_m$ and the relation
$a\phimeq x$ holds for the only $x=\phi_{m+1}(a)$.
We have $a'\ne\phi_i(a)=\phi_{m+1}(a)$ (the latter equality is due to
the uniqueness of the $\phi_i$'s ensured by Condition 1).
Therefore $a\notphimeq a'$, which means that $(a,a')$ is $m$-threatening.
\end{subproof}

\begin{claim}\label{cl:threat}
Assume that Phase $j$, $j\le k+1$, finishes, Conditions 1--3 for all
$i\le j$ are met, and the game is going on. Let $1\le i\le j$.
As soon as after Phase $j$ an $i$-threatening pair $(a,a')$ arises,
Spoiler is able to win in at most $(i-1)(k-1)$ moves playing all the time,
at his own choice, either in $M$ or in $M'$.
\end{claim}

\noindent
{\sl Convention.} Given a relation $R=R^M$ of $M$, we will denote
the respective relation $R^{M'}$ by $R'$.

\medskip

\begin{subproof}
We proceed by induction on $i$. For $i=1$ the claim easily follows from
Item 3 of Lemma \ref{lem:phieq}. Let $i\ge2$ and assume that the claim
is true for all preceding values $1,2,\ldots,i-1$.

We focus on the case that $a\notin X_i$ (in the case that $a'\notin X'_i$
the proof is given by the symmetric argument). The non-equivalence
$a\notphiieq a'$ can happen in two situations.

\case 1{$a'\in X'_i$.}
Clearly, $a\ne\phi_i^{-1}(a')$ and therefore, by Claim \ref{cl:conserv},
the pair $(a,a')$ is $m$-threatening for some $m<i$. By the induction
hypothesis, Spoiler is able to win in at most $(m-1)(k-1)$ moves.

\case 2{$a'\notin X'_i$.}
Then the non-equivalence $a\notphiieq a'$ means that there is an
$l$-ary relation of $M$ and $\bar a\in(X_i\cup\{a\})^l$ with at least
one occurrence of $a$ such that
\begin{equation}\label{eq:psibad}
R\bar a\ne R'\psi\bar a,
\end{equation}
where $\psi$ is the map defined by $\psi(x)=\phi_i(x)$ for all
$x$ in $A=\{a_1,\ldots,a_l\}\setminus\{a\}$ and by $\psi(a)=a'$.
Thus, $\psi$ is not a partial isomorphism from $M$ to $M'$.
Hence, if $A\subseteq\dom\phi^*_i$, then Spoiler wins immediately.

Assume that $\hat A=A\setminus\dom\phi^*_i$ is nonempty.
Spoiler selects all unselected elements in $\hat A$, if he wants to
play in $M$, or in $\phi_i(\hat A)$, if he prefers to play in $M'$.
This takes at most $k-1$ moves. Suppose that Spoiler plays in $M$
(for $M'$ the argument is symmetric).
If for every $b\in\hat A$ its counterpart in $V(M')$
is $\phi_i(b)$, this is Spoiler's win by \refeq{eq:psibad}.
If some $b\in\hat A$ has the counterpart $b'$ such that $b'\ne\phi_i(b)$,
by Claim \ref{cl:conserv} there arises an
$m$-threatening pair for some $m<i$.
Applying the induction hypothesis for the index $m$, we conclude that
Spoiler is able to win in at most $(m-1)(k-1)$ moves, having made altogether
at most $(k-1)+(m-1)(k-1)\le(i-1)(k-1)$ moves.
\end{subproof}

\begin{claim}\label{cl:cond}
Assume that Phase $j$, $j\le k$, has been finished and
Conditions 1--3 for all $i\le j$ are met.
Assume furthermore that Part 1 of Phase $j+1$ finishes and the game
is still going on.
Then either Conditions 1--3 hold true for $i=j+1$ as well or Spoiler
is able to win or to create an $i$-threatening pair for some $i\le j$
in at most $k$ moves with at most one alternation from $M$ to $M'$
(and hence he is able to win in Part 2 of Phase $j+1$).
\end{claim}

\begin{subproof}
Assuming that Spoiler is unable to win or to create an $i$-threatening pair,
we check Conditions 1--3.

{\em Condition 2.}
The following two facts take place, for else
Spoiler would be able to enforce creating a $j$-threatening pair
in at most one move:
\begin{itemize}
\item
For every $C'\in\calC^{k+1}(X'_j)$ there is $C\in\calC^{k+1}(X_j)$
such that $C\phijeq C'$. (Otherwise Spoiler selects an elements in $C'$
that violates this condition and a $j$-threatening pair arises
whatever Duplicator's response.)
\item
For every $C\in\calC^{k+1}(X_j)$ there is $C'\in\calC^{k+1}(X'_j)$
such that $C\phijeq C'$ and $|C'|\ge|C|-1$.
(Otherwise, for some $C$, $\phi^*_{j+1}(\tilde C)$ cannot be included
into the respective $C'$. Therefore $c\notphijeq \phi^*_{j+1}(c)$ for
at least one $c\in\tilde C$, providing us with a $j$-threatening pair.)
\end{itemize}
Thus, there is a one-to-one correspondence between $\calC^{k+1}(X_j)$
and $\calC^{k+1}(X'_j)$ such that, for $C$ and $C'$ corresponding to one
another, $C\phijeq C'$, $|C'|\ge|C|-1$, and
$\phi^*_{j+1}(\tilde C)\subseteq C'$. Moreover, it actually holds $|C'|=|C|$
because, if $|C'|\ge|C|+1$, Spoiler could select 2 vertices in
$C'\setminus\phi^*_{j+1}(\tilde C)$ obtaining a $j$-threatening
pair whatever Duplicator's response.

{\em Conditions 1 and 3.}
By Condition 1 for $i=j$, the partial isomorphism $\phi^*_{j+1}$ can
be extended on $X_j$ only to $\phi_j$ and then it remains undefined
within $X_{j+1}$ only on $Y_j\setminus\tilde Y_j$. Define an extension
$\phi_{j+1}$ of $\phi^*_{j+1}$ on the whole $X_{j+1}$ so that
$\phi_{j+1}$
\begin{itemize}
\item
agrees with $\phi_j$ on $X_j$,
\item
agrees with $\phi^*_{j+1}$ on $\tilde Y_j$, and
\item
for each $C\in\calC^{k+1}(X_j)$, takes the single element in
$C\setminus\tilde C$ to
the single element in $C'\setminus\phi^*_{j+1}(\tilde C)$, where $C'$
corresponds to $C$ according to Condition 2 that we have already proved.
\end{itemize}
We have to show that $\phi_{j+1}$ is a partial isomorphism from $M$
to $M'$ and no other extension of $\phi^*_{j+1}$ is such.

Assume that $\phi_{j+1}$ is not a partial isomorphism and get a
contradiction to the assumption that Spoiler can in the nearest $k$
moves neither win nor create an $i$-threatening pair.
For some $l$-ary relation of $M$ and $\bar a\in X^l_{j+1}$,
we should have
\begin{equation}\label{eq:notiso}
R\bar a\ne R'\phi_{j+1}\bar a.
\end{equation}
As a consequence, $A=\{a_1,\ldots,a_l\}$ is not included into
$\dom\phi^*_{j+1}$ for else $\phi^*_{j+1}$ would not be a partial
isomorphism, contradicting the assumption that the game is still
going on. Let Spoiler select all elements in
$\hat A=A\setminus\dom\phi^*_{j+1}$. If for $b\in\hat A$ Duplicator
always responds with $\phi_{j+1}(b)$, he loses by \refeq{eq:notiso}.
Otherwise, let $b$ be an element in $\hat A$ to which Duplicator
responds with $b'\ne\phi_{j+1}(b)$. If $b\in X_j$ or $b'\in X'_j$, then
we have $b'\ne\phi_j(b)$ because $\phi_{j+1}$ extends $\phi_j$, .
By Claim \ref{cl:conserv}, $(b,b')$ is an $i$-threatening pair for
some $i<j$. If $b\in X_{j+1}\setminus X_j$ and
$b'\in X'_{j+1}\setminus X'_j$, then $b\notphijeq b'$ by Condition 2
proved above and the definition of $\phi_{j+1}$. Thus, $(b,b')$ is
$j$-threatening. We have a contradiction in any case and therefore
$\phi_{j+1}$ is a partial isomorphism from $M$ to $M'$ indeed.

To prove the uniqueness of the extension $\phi_{j+1}$ (i.e., Condition 1),
assume that $\hat\phi_{j+1}$ is another extension of $\phi^*_{j+1}$
over $X_{j+1}$
which is a partial isomorphism and differs from $\phi_{j+1}$
at $b\in Y_j\setminus\tilde Y_j$. Let $b'=\phi_{j+1}(b)$ and
$b''=\hat\phi_{j+1}(b)$. By Condition 2 proved above,
\begin{equation}\label{eq:noteqxjj}
b'\notxjjeq b''.
\end{equation}
By Condition 1 for $i=j$, $\hat\phi_{j+1}$ on $X_j$ coincides with
$\phi_j$. Thus, the composition $\hat\phi_{j+1}\phi^{-1}_{j+1}$
takes $b'$ to $b''$, extends $\id_{X'_j}$, and is an automorphism
of $M'[X'_{j+1}]$. This makes a contradiction to~\refeq{eq:noteqxjj}.
\end{subproof}

Claim \ref{cl:cond} implies by an easy induction on $j$ from 1 to $k+1$
that, for each $1\le j\le k$, unless Spoiler wins in Phase $j$ or
earlier, Conditions 1--3 assumed in our description of Phase $j+1$
are indeed true. For analysis of the concluding phase, we state
simple consequences of Claims~\ref{cl:conserv}--\ref{cl:cond}.

\begin{claim}\label{cl:concl}
Suppose that Spoiler follows the strategy designed above
(Duplicator's strategy does not matter). Assume that Duplicator
survives up to Phase $k+1$. Then the following claims are true.
\begin{enumerate}
\item
Conditions 1--3 hold true for all $i\le k+1$.
\item
When in further play Spoiler selects $v\in V(M)\cup V(M')$,
we denote Duplicator's response by $\psi(v)$.
As long as there arises no $i$-threatening pair for any $i\le k$,
it holds
\begin{eqnarray}
\psi(v)\phikeq v&\mbox{\ if\ }&v\notin X_k\cup X'_k,\label{eq:notin}\\
\psi(v)=\phi_{k+1}(v)&\mbox{\ if\ }&v\in X_{k+1},\label{eq:xm}\\
\psi(v)=\phi_{k+1}^{-1}(v)&\mbox{\ if\ }&v\in X'_{k+1}.\label{eq:xmm}
\end{eqnarray}
(The relations in \refeq{eq:notin} and \refeq{eq:xm}--\refeq{eq:xmm}
are equivalent on $(X_{k+1}\cup X'_{k+1})\setminus
(X_k\cup X'_k)$.)
\end{enumerate}
\end{claim}

\begin{subproof}
Item 1 follows from Claim \ref{cl:cond} by an easy inductive argument.
Regarding Item 2, note that, if \refeq{eq:notin} were false, $(v,\psi(v))$
would be a $k$-threatening pair. If \refeq{eq:xm} or \refeq{eq:xmm}
were false, $(v,\psi(v))$ would be an $i$-threatening pair for some
$i\le k$ on the account of Claim~\ref{cl:conserv}.
\end{subproof}

\bigskip

{\sc Concluding Phase (Phase $k+2$).}

\smallskip

We here assume that Phases from 1 up to $k+1$ have been finished
without Spoiler's win and therefore Items 1 and 2 of Claim \ref{cl:concl}
hold true. As soon as there arises an $i$-threatening pair
for some $i\le k$, Spoiler switches to the strategy given by
Claim \ref{cl:threat} and wins in at most $(k-1)^2$ moves.
As long as there occurs no such pair, Spoiler follows the strategy
described below. The strategy depends on which of the following
three cases takes place.

\smallskip

\case 1{{\em
There is a one-to-one correspondence between $\calC(X_k)$ and $\calC(X'_k)$
such that, if $C$ and $C'$ correspond to one another, then
$C\phikeq C'$ and, moreover, $|C|=|C'|$.
}}
By Item 1 of Claim \ref{cl:concl}, such correspondence does exist
between $\calC^{k+1}(X_k)$ and $\calC^{k+1}(X'_k)$ in any case.

Let $\Upsilon$ be the set of maps $\phi\function{V(M')}{V(M)}$ such that
\begin{itemize}
\item
$\phi$ is one-to-one,
\item
$\phi$ extends $\phi_{k+1}^{-1}$,
\item
for every $C'\in\calC(X'_k)$, we have $\phi(C')\in\calC(X_k)$ and
$\phi(C')\phikeq C'$.
\end{itemize}

\begin{claim}\label{cl:upsilon}
Assume that $\phi$ and $\psi$ are in $\Upsilon$. Let $R$ be an $l$-ary
relation of $M$. Then $R\phi\bar a'=R\psi\bar a'$ for all
$\bar a'\in V(M')^l$.
\end{claim}

\begin{subproof}
The product $\psi\phi^{-1}$ is a permutation of $V(M)$ that moves
only elements in $Z$. Moreover, $\psi\phi^{-1}$ preserves the partition
$\calC(X_k)\setminus\calC^{k+1}(X_k)$ of $Z$ and therefore $\psi\phi^{-1}$
is decomposed into the product of permutations $\pi_C$ over
$C\in\calC(X_k)\setminus\calC^{k+1}(X_k)$, where each $\pi_C$ acts
on the respective $C$. Since every $\pi_C$ is decomposable into
a product of transpositions, we have $\psi\phi^{-1}=\tau_1\tau_2\ldots\tau_t$
with $\tau_i$ being a transposition of two elements both in some $C$.
It is easy to see that
$\psi\bar a'=(\ldots((\phi\bar a')^{\tau_t})\ldots)^{\tau_1}$.
By Lemma \ref{lem:xeqsim}, each application of $\tau_i$ does not change
the initial value of $R\phi\bar a'$. Therewith we arrive at the desired
equality $R\phi\bar a'=R\psi\bar a'$.
\end{subproof}

To specify Spoiler's strategy, we fix $\phi\in\Upsilon$ arbitrarily.
Since $M$ and $M'$ are nonisomorphic, $\phi$ is not an isomorphism
from $M'$ to $M$, that is,
\begin{equation}\label{eq:phinotiso}
R\phi\bar a'\ne R'\bar a'
\end{equation}
for some $l$-ary relation $R'$ of $M'$
and $\bar a'\in V(M')^l$.
This inequality implies that the set $A'=\{a'_1,\ldots,a'_l\}$
is not included into $X'_{k+1}$. Spoiler selects, one by one,
elements of $\hat A=A'\setminus\range\phi^*_{k+1}$. For Spoiler's
move $v$, let $\psi(v)$ denote Duplicator's response.

Assume first that
\begin{eqnarray}
\psi(v)\phikeq v&\mbox{\ whenever\ }&v\notin X'_k\nonumber\\
&\mbox{and}&\label{eq:psip}\\
\psi(v)=\phi_{k+1}^{-1}(v)&\mbox{\ whenever\ }&v\in X'_{k+1}.\nonumber
\end{eqnarray}
Due to \refeq{eq:psip}, we are able to
extend $\psi$, initially defined on $\hat A$, to a map in $\Upsilon$.
Fix a such extension. By Claim \ref{cl:upsilon},
$R\phi\bar a'=R\psi\bar a'$ and, by \refeq{eq:phinotiso},
Spoiler wins. If \refeq{eq:psip} is violated for some $v\in\hat A$,
by Item 2 of Claim \ref{cl:concl} this produces an $i$-threatening
pair for $i\le k$ and therefore Spoiler wins in at most $(k-1)^2$ moves,
having made altogether at most $k+(k-1)^2$ moves.

\case 2{{\em
There is no one-to-one correspondence between $\calC(X_k)$ and $\calC(X'_k)$
such that, if $C$ and $C'$ correspond to one another, then $C\phikeq C'$.
}}
Spoiler selects an element in $C$ or $C'$ that has no counterpart.
If Duplicator responds with a vertex outside $X_k\cup X'_k$,
there arises a $k$-threatening pair.
If Duplicator responds with a vertex in $X_k\cup X'_k$,
there arises an $i$-threatening pair for $i<k$ by Item 2 of Claim
\ref{cl:concl}.
This allows Spoiler to win altogether in at most $1+(k-1)^2$ moves.

\case 3{{\em
There is a one-to-one correspondence between $\calC(X_k)$ and $\calC(X'_k)$
such that, if $C$ and $C'$ correspond to one another, then $C\phikeq C'$.
However, there are $C\in\calC(X_k)$ and $C'\in\calC(X'_k)$ such that
$C\phikeq C'$ but $|C|\ne|C'|$.
}}

Call a class $C\in\calC(X_k)$ {\em useful\/}
if $C\phikeq C'$ but $|C|\ne|C'|$. The description
of Case 3 tells us that there is at least one useful class. Actually,
since $|V(M)|=|V(M')|$, there are at least two useful classes,
$C_1$ and $C_2$. Note that $|C_1|+|C_2|\le|Z|$. Without loss of generality,
assume that $|C_1|\le|Z|/2$. Let $C'_1$ be the counterpart of $C_1$
in $\calC(X'_k)$, i.e., $C_1\phikeq C'_1$.
In the larger of $C_1$ and $C'_1$ Spoiler selects $\min\{|C_1|,|C'_1|\}+1$
elements.
Duplicator is enforced to at least once reply not in the smaller class.
By Item 2 of Claim \ref{cl:concl}, this produces an $i$-threatening pair
and Spoiler, according to Claim \ref{cl:threat}, wins in at most
$(k-1)^2$ subsequent moves, having made altogether at most
$|Z|/2+1+(k-1)^2$ moves.

{\sc End of description of the concluding phase}

\subsection{Estimation of the length of the game}\label{ss:length}

If Spoiler follows the above strategy and Duplicator delays his loss
as long as possible, the end of the game is always this:
Spoiler enforces creating a threatening pair in at most $k$ moves
and then wins in at most $(k-1)^2$ next moves using the strategy
of Claim \ref{cl:threat}. Let us calculate the smallest possible
(optimal for Duplicator) number of elements in $M$ unoccupied
till such final stage of the game. The minimum is attained if
all Phases from 1 up to $k+2$ are played and in Phase $k+2$
it happens Case 3. Then the number of elements unoccupied in $X_{k+1}$
is equal to
$$
\sum_{i=1}^k |Y_i\setminus\tilde Y_i|=\sum_{i=1}^k |\calC^{k+1}(X_i)|.
$$
The number of elements unoccupied in $Z$ is at least
$|Z|-(|Z|/2+1)=|Z|/2-1$.
By Lemma \ref{lem:unpebbl}, the total number of unoccupied elements
is at least
$$
\sum_{i=1}^k |\calC^{k+1}(X_i)|+\frac{|Z|}2-1 >
\frac n{2k}-\frac12-\frac1{2k}.
$$
Thus, the maximum possible number of occupied elements is less than
$$
\of{1-\frac1{2k}}n+\frac12+\frac1{2k}.
$$
Summing up, we conclude that our strategy allows Spoiler to win
in less that
$$
\of{1-\frac1{2k}}n+k^2-k+2
$$
moves. Theorem \ref{thm:main} is proved.

\subsection{Definability results}

A natural question is if our approach applies to defining rather than
identifying formulas. In fact,
the proof of Theorem \ref{thm:main} implies the definability with
lower quantifier rank for a quite representative class of structures.

\subsubsection{Definability of irredundant structures}

\begin{definition}\label{def:sigma}
If $M$ is a finite structure, let
$$
\gs(M)=\max\setdef{|A|}{A\subseteq V(M)\mbox{\ such\ that\ }
a_1\sim a_2\mbox{\ for\ every\ }a_1,a_2\in A}
$$
be the maximum cardinality of a $\sim$-equivalence class in~$V(M)$.

If $\gs(M)=1$, i.e., no transposition of two elements is an automorphism
of $M$, we call $M$ {\em irredundant}.
\end{definition}

\begin{theorem}\label{thm:irred}
Let $M$ be an irredundant structure of order $n$ with maximum relation
arity $k$. Then
$$
\Dfn 1M < \of{1-\frac1{2k}}n+k^2-k+1.
$$
\end{theorem}

\begin{proof}
It is not hard to see that an irredundant structure whose all relations
are unary is definable by a formula with quantifier rank 1.
Assume therefore that $k\ge2$.
Notice that Spoiler's strategy described in Section \ref{ss:strategy}
applies for any pair of $L$-structures $M$ and $M'$ of arbitrary orders
with the only exception of Case 3 in the concluding Phase $k+2$,
where the equality $|V(M)|=|V(M')|$ is supposed. Since the set $Z$
is partitioned into $\sim$-equivalence classes each consisting of
at least $k+2$ elements, for an irredundant structure $M$ we have
$Z=\emptyset$. Consequently, $V(M)=X_{k+1}$. It follows that either
Spoiler wins at latest in Phase $k+1$ or, according to Item 1 of
Claim \ref{cl:concl}, there is a partial isomorphism $\phi_{k+1}$
from $M$ to $M'$ with $\dom\phi_{k+1}=V(M)$.

In the latter case,
since $M$ and $M'$ are non-isomorphic, there is at least one element
$v\in V(M')\setminus\range\phi_{k+1}$. In the concluding phase of the
game Spoiler selects $v$
and, according to Item 2 of Claim \ref{cl:concl}, there arises
a $k$-threatening pair. Spoiler  switches to the strategy given by
Claim \ref{cl:threat} and wins in at most $(k-1)^2$ moves.

It remains to estimate the length of the game. Similarly to
Section \ref{ss:length}, we conclude that Spoiler needs at most
$n-\sum_{i=1}^k|\calC^{k+1}(X_i)|+k+(k-1)^2$ to win. By estimate
\refeq{eq:lastbut}, where $|Z|=0$, this number is less than
$(1-\frac1{2k})n+k^2-k+1$.
\end{proof}

\begin{remark}\label{rem:irred}
There are simple examples of irredundant structures $M$ showing
a lower bound $\dfn M\ge n/4$. For example, let $F$ be a directed
graph on two vertices $u$ and $v$ consisting of a single (directed)
edge $(uv)$. Let $G$ be another directed graph on $u$ and $v$
consisting of two edges, $(uv)$ and the loop $(uu)$. Denote the
disjoint union of $a$ copies of $F$ and $b$ copies of $G$ by $aF+bG$.
It is easy to see that $aF+bG$ is irredundant for any $a$ and $b$.
Directed graphs $M=mF+mG$ and $M'=(m-1)F+(m+1)G$ are non-isomorphic
and both have order $4m$. An obvious strategy for Duplicator in the
\EF\/ game on $M$ and $M'$ shows that $\D{M,M'}\ge m$.
\end{remark}

Theorem \ref{thm:irred} will be considerably strengthened in the
next subsections. In particular, it will be surpassed by
Theorem~\ref{thm:dmupper}.

\subsubsection{A further refinement}

As we observed in the proof of Theorem \ref{thm:irred}, Spoiler's
strategy designed in Section \ref{ss:strategy} ensures the bound
\begin{equation}\label{eq:dmm}
\DD1{M,M'} < \of{1-\frac1{2k}}n+k^2-k+2.
\end{equation}
for $M'$ of {\em any\/} order under an additional condition imposed on $M$.
We are able to describe exceptional pairs of non-isomorphic $M$ and $M'$
for which \refeq{eq:dmm} may not hold much more precisely.
Assume that $M'$ has order $n'\ge n$. As was already mentioned,
the assumption that $n'=n$ is used only in Case 3 of the concluding
Phase $k+2$. Turning back to this case, we see that what is actually
used is the existence of at least two useful classes in $\calC(X_k)$.
Thus, \refeq{eq:dmm} may not hold in the only case that there is
a unique useful class $C_0\in\calC(X_k)$. Since actually
$C_0\in\calC(X_k)\setminus\calC^{k+1}(X_k)$, we have $|C_0|\ge k+2$.
By Lemma \ref{lem:xeqsim}, the class $C_0$ consists of pairwise
$\sim$-equivalent elements.

Let $C'_0$ be the counterpart of $C_0$ in $\calC(X'_k)$, i.e.,
$C'_0\phikeq C_0$. Given $B\subseteq C'_0$ with $|B|=|C_0|$,
let $M'_B=M'[V(M')\setminus(C'_0\setminus B)]$. Consider an arbitrary map
$\phi\function{V(M'_B)}{V(M)}$ extending $\phi^{-1}_{k+1}$,
mapping each $C'\in\calC(X'_k)\setminus\{C'_0\}$
onto its $\phikeq$-counterpart in $\calC(X_k)$, and
mapping $B$ onto $C_0$.
As in Case 1 of Phase $k+2$, we see that Spoiler is able to win
within the bound of \refeq{eq:dmm} unless $\phi$ is an isomorphism
from $M'_B$ to $M$. From here we easily arrive at the following
conclusion.

\begin{lemma}\label{lem:dmmpiso}
Let $L$ be a vocabulary with maximum relation arity $k$.
Let $M$ and $M'$ be non-isomorphic $L$-structures of orders
$n$ and $n'$ respectively and $n\le n'$. Then the bound
$$
\DD1{M,M'} < \of{1-\frac1{2k}}n+k^2-k+2
$$
may be false only if there is a set $C_0\subseteq V(M)$ with
$|C_0|\ge k+2$ consisting of pairwise $\sim$-equivalent vertices
and there is a partial isomorphism $\psi$ from $M$ to $M'$
defined on $V(M)\setminus C_0$ whose any injective extension
is a partial isomorphism from $M$ to~$M'$.
\end{lemma}

In the next subsection we make a constructive interpretation of
the condition appearing in the lemma.

\subsubsection{Cloning an element of a structure}\label{sss:clone}

{\bf Notation.}
Recall that, given a set $V$ and a function $\pi$ defined on $V$,
we extend $\pi$ over $V^l$, where $l\ge1$, by
$\pi\baru=(\pi(u_1),\ldots,\pi(u_l))$ for any $\baru=(u_1,\ldots,u_l)$
with all $u_i$ in $V$. In particular, this concerns the case that
$\pi$ is a permutation of elements of $V$. Recall also that, if
$\pi=(v_1v_2)$ is a transposition, then we may write $\baru^{(v_1v_2)}$
in place of $\pi\baru$.

\begin{definition}
Given $v\in V(M)$, let $[v]_M=\setdef{u\in V(M)}{u\sim v}$ be the
$\sim$-equiva\-le\-n\-ce class of the element $v$.
\end{definition}

We now introduce an operation of expanding a class $[v]_M$, i.e.,
adding to $M$ new elements $\sim$-equivalent to $v$.
This operation was considered in \cite{PVV} in the particular case
of uniform hypergraphs.

Let $L$ be a vocabulary with maximum relation arity $k$.
Below $K$ and $M$ are $L$-structures, $v$ is an element of $M$,
and $t$ is a non-negative integer.

\medskip

\noindent
{\bf Definition A}
The notation $K=M\add tv$ means that the following conditions are
fulfilled.
\begin{description}
\item[A1]
$V(M)\subseteq V(K)$ and $|V(K)|=|V(M)|+t$.
\item[A2]
$K[V(M)]=M$.
\item[A3]
$|[v]_M|\ge k$.
\item[A4]
$[v]_K=[v]_M\cup(V(K)\setminus V(M))$.
\end{description}

\medskip

\noindent
{\bf Definition B}
The notation $K=M\add tv$ means that the following conditions are
fulfilled.
\begin{description}
\item[B1]
$V(M)\subseteq V(K)$ and $|V(K)|=|V(M)|+t$.
\item[B2]
There is $C\subseteq[v]_M$ with $|C|\ge k$ such that every injective
extension of $\id_{V(M)\setminus C}$ to a map
$\psi\function{V(M)}{V(K)}$ is a partial isomorphism from $M$ to~$K$.
\end{description}

\medskip

\noindent
{\bf Definition C}
The notation $K=M\add tv$ means that the following conditions are
fulfilled.
\begin{description}
\item[C1]
$V(M)\subseteq V(K)$ and $|V(K)|=|V(M)|+t$.
\item[C2]
$|[v]_M|\ge k$.
\item[C3]
Let $R$ be an $l$-ary relation in $L$. If $\baru\in V(M)^l$, then
$R^K\baru=R^M\baru$.
\item[C4]
Let $R$ be an $l$-ary relation in $L$.
Assume that $\baru\in V(K)^l$ and the set
$\{u_1,\ldots,u_l\}\setminus V(M)=\{w_1,\ldots,w_p\}$ is nonempty.
Then $R^K\baru=1$ iff there are pairwise distinct elements
$v_1,\ldots,v_p\in[v]_M\setminus\{u_1,\ldots,u_l\}$ such that
$R^M\pi\baru=1$ for $\pi=(w_1v_1)\cdots(w_pv_p)$.
\end{description}

\begin{lemma}\label{lem:defequi}
Definitions A, B, and C are equivalent.
\end{lemma}

\begin{proof}
{\em Conditions A1--A4 imply Conditions B1--B2}.
Since B1 coincides with A1, we only have to derive B2. We are
actually able to prove B2 for an arbitrary $C\subseteq[v]_M$
with $|C|\ge k$ (there is at least one such $C$ by A3). Let
$\psi$ be as specified in B2. For any $l$-ary relation $R$ in $L$
and $\baru\in V(M)^l$, we have to check that $R^M\baru=R^K\psi\baru$.
Assume that in $\{\psi(u_1),\ldots,\psi(u_l)\}$ there are $p$
elements from $V(K)\setminus V(M)$ and denote them by $w_1,\ldots,w_p$.
Take arbitrary pairwise distinct
$v_1,\ldots,v_p\in C\setminus\{\psi(u_1),\ldots,\psi(u_l)\}$.
Let $\tiu=\pi\psi\baru$ with $\pi=(w_1v_1)\cdots(w_pv_p)$.
By A4, we have $v_i\sim w_i$ in $K$ for all $i\le p$.
It follows that $R^K\psi\baru=R^K\tiu$. Since $\tiu\in V(M)^l$,
by A2 we have $R^K\tiu=R^M\tiu$. Notice now that $\tiu$ and $\baru$
coincide at the positions occupied by elements in $V(M)\setminus C$,
while elements in $C$ are permuted according to some permutation
$\tau$, i.e., $\tiu=\tau\baru$. Since $\tau$ is decomposable in a
product of transpositions and elements of $C$ are pairwise $\sim$-equivalent
in $M$, we have $R^M\tiu=R^M\baru$, completing derivation of~B2.

\medskip

{\em Conditions B1--B2 imply Conditions C1--C4}.
For C1 and C2 this is trivial.
C3 immediately follows from B2 if we take $\psi=\id_{V(M)}$.
Let us focus on C4. Let $\baru$ and $w_1,\ldots,w_p$ be as specified
in this condition. Assume first that $R^K\baru=1$. Take
$v_1,\ldots,v_p\in C\setminus\{u_1,\ldots,u_l\}$ being pairwise
distinct and define $\psi$ by $\psi(v_i)=w_i$ for $i\le p$
and $\psi(x)=x$ for all other $x\in V(M)$. Notice that
$\psi^{-1}\baru=\pi\baru$ for $\pi=(w_1v_1)\cdots(w_pv_p)$.
As $\psi$ is a partial isomorphism by B2, we conclude that
$R^M\pi\baru=R^K\baru=1$.
This proves C4 in one direction. Such a way of proving
$R^M\pi\baru=R^K\baru$ will be referred to as {\em $\psi$-argument}.

For the other direction, assume that $R^M\pi\baru=1$ for
$\pi=(w_1v_1)\cdots(w_pv_p)$ with some
$v_1,\ldots,v_p\in [v]_M\setminus\{u_1,\ldots,u_l\}$.
If all $v_i$ are in $C$, the equality $R^K\baru=1$ follows from the
$\psi$-argument with the same $\psi$ as above.
Otherwise, we can replace each $v_i$ with some
$v'_i\in C$, where $v'_1,\ldots,v'_p$ are pairwise distinct
elements of $C\setminus\{u_1,\ldots,u_l\}$ and $v'_i=v_i$
whenever $v_i\in C$. For no $i$ this
replacement changes the initial value of $R^M\pi\baru$ and,
after all replacements are done, we have $R^M\pi'\baru=1$
with $\pi'=(w_1v'_1)\cdots(w_pv'_p)$.
Defining $\psi'$ by $\psi'(v'_i)=w_i$ and $\psi'(x)=x$ elsewhere
on $V(M)$, we obtain $R^K\baru=R^M\pi'\baru=1$ by the $\psi'$-argument.

\medskip

{\em Conditions C1--C4 imply Conditions A1--A4}.
Since A1--A3 are virtually the same as C1--C3, our concern is A4.
It is easy to see
that $[v]_K\cap V(M)$ cannot be larger than $[v]_M$. Therefore,
it suffices to show that in $K$ we have $v\sim v'$ for any
$v'\in [v]_M\cup(V(K)\setminus V(M))$. Given an $l$-ary relation
$R$ in $L$ and $\baru\in V(K)^l$, we have to check that
$$
R^K\baru=R^K\baru^{(vv')}.
$$
We do it by routine examination of several cases. Note that, if neither
$v$ nor $v'$ occurs in $\baru$, then there is nothing to prove.

To simplify notation, denote
$$\hatu=\baru^{(vv')}.
$$
Furthermore, let $U=\{u_1,\ldots,u_l\}$ and
$U\setminus V(M)=\{w_1,\ldots,w_p\}$.
Denote the set of elements in $\hatu$ by $\hatU$.

\case 1{$v'\in V(K)\setminus V(M)$.}

\subcase{1.1}{$v\in U$, $v'\in U$.}\\
Assuming $R^K\baru=1$, we will infer $R^K\hatu=1$. This will give also
the converse implication because $\baru$ is supposed arbitrary
with occurrences of both $v$ and $v'$ and we hence can take $\hatu$ instead
of $\baru$. Without loss of generality, assume that $v'=w_p$. By C4,
there are $v_1,\ldots,v_p\in[v]_M\setminus U$ such that
$R^M\pi\baru$ with $\pi=(w_1v_1)\cdots(w_{p-1}v_{p-1})(v'v_p)$.
As easily seen, $\pi\hatu=(\pi\baru)^{(vv_p)}$.
Since $v_p\sim v$ in $M$, we have $R^M\pi\hatu=1$.
Note that $\hatU=U$ and hence $v_1,\ldots,v_p\in[v]_M\setminus\hatU$.
By C4, we conclude that $R^K\hatu=1$, as desired.

\subcase{1.2}{$v\in U$, $v'\notin U$.}\\
Note that $\hatU\setminus V(M)=\{w_1,\ldots,w_p,v'\}$ and
$[v]_M\setminus\hatU=([v]_M\setminus U)\cup\{v\}$.
We first assume that $R^K\baru=1$ and infer from here that
$R^K\hatu=1$. Let $v_1,\ldots,v_p$ be as ensured by Condition C4
for $\baru$, that is, $R^M\pi\baru=1$ with $\pi=(w_1v_1)\cdots(w_pv_p)$.
Let $\pi'=\pi(v'v)$. As easily seen, $\pi'\hatu=\pi\baru$.
Thus, $R^M\pi'\hatu=R^M\pi\baru=1$ and, by C4, we conclude that
$R^K\hatu=1$.

We now assume that $R^K\hatu=1$ and have to infer $R^K\baru=1$.
According to C4, there are pairwise distinct
$v'_1,\ldots,v'_{p+1}\in[v]_M\setminus\hatU$ such that
$R^M\pi'\hatu=1$ with $\pi'=(w_1v'_1)\cdots(w_pv'_p)(v'v'_{p+1})$.
Choose pairwise distinct $v_1,\ldots,v_p$ in
$\{v'_1,\ldots,v'_{p+1}\}\setminus\{v\}$ and apply to $\baru$
the substitution $\pi=(w_1v_1)\cdots(w_pv_p)$. It is not hard to see
that $\pi\baru=\tau\pi'\hatu$ for $\tau$ being a permutation of
the set $V=\{v,v'_1,\ldots,v'_p,v'_{p+1}\}$ taking $v'_i$ to $v_i$
for $i\le p$ and $v'_{p+1}$ to $v$. A such $\tau$ exists because
elements in $\{v'_1,\ldots,v'_{p+1}\}$ and in $\{v_1,\ldots,v_p,v\}$
are pairwise distinct (the fact that the two sets may intersect
does not matter). Since $\tau$ is decomposable in a product of
transpositions of two elements from $V$ and elements in $V$ are
pairwise $\sim$-equivalent in $M$, we have
$R^M\pi\baru=R^M\pi'\hatu=1$. By C4, we conclude that
$R^K\baru=1$, as desired.

\subcase{1.3}{$v\notin U$, $v'\in U$.}\\
This subcase reduces to Subcase 1.2 by considering $\hatu$
in place of~$\baru$.

\case 2{$v'\in[v]_M$.}\\
Since in this case $v$ and $v'$ are interchangeable, it suffices to
assume that $v\in U$ and prove that $R^K\baru=1$ implies $R^K\hatu=1$.
Note that $\hatU\setminus V(M)=\{w_1,\ldots,w_p\}$.

\subcase{2.1}{$v'\in U$.}\\
Note that $[v]_M\setminus\hatU=[v]_M\setminus U$.
Let $v_1,\ldots,v_p\in[v]_M\setminus U$ be as ensured by Condition C4
for $\baru$, i.e., $R^M\pi\baru=1$ with $\pi=(w_1v_1)\cdots(w_pv_p)$.
Applying the same $\pi$ to $\hatu$, we see that
$\pi\hatu=(\pi\baru)^{(vv')}$. As $v\sim v'$ in $M$, we have
$R^M\pi\hatu=R^M\pi\baru=1$ and hence, by C4, we obtain $R^K\hatu=1$.

\subcase{2.2}{$v'\notin U$.}\\
Note that $[v]_M\setminus\hatU=(([v]_M\setminus U)\setminus\{v'\})\cup\{v\}$.
Let $v_1,\ldots,v_p$ and $\pi$ be as in Subcase 2.1. The difference is
that now the containment $v'\in\{v_1,\ldots,v_p\}$ is possible.
For $i\le p$, set
$$
v'_i=\cases{
v_i & if $v_i\ne v'$\cr
v & if $v_i=v'$\cr}
$$
and apply to $\hatu$ the substitution $\pi'=(w_1v'_1)\cdots(w_pv'_p)$.
It is not hard to see that $\pi'\hatu=\tau\pi\baru$ for
$\tau$ being a permutation of the set $\{v,v_1,\ldots,v_p,v'\}$
taking $v_i$ to $v'_i$ for all $i\le p$ and $v$ to $v'$.
Similarly to the second part of Subcase 1.2, we conclude that
$R^M\pi'\hatu=R^M\pi\baru=1$ and, by C4, we obtain $R^K\hatu=1$.
\end{proof}

\begin{lemma}\label{lem:clone}
Let $L$ be a vocabulary with maximum relation arity $k$.
Let $M$ be an $L$-structure, $v\in V(M)$ with $|[v]_M|\ge k$,
and $t\ge0$. Then an $L$-structure $K$ such that $K=M\add tv$
exists and is unique up to an isomorphism.
\end{lemma}

\begin{proof}
The existence follows from Definition C. To obtain $K$, we add $t$
new elements to $V(M)$, keep all relations of $M$ on $V(M)$, and
add new relations involving at least one new element, being
guided by Condition~C4.

To prove the uniqueness, we use Definition B. Assume that
$K_1=M\add tv$ and $K_2=M\add tv$ according to this definition.
Let $\phi\function{V(K_1)}{V(K_2)}$ be an arbitrary one-to-one
map whose restriction on $V(M)$ is $\id_{V(M)}$.
We claim that $\phi$ is an isomorphism from $K_1$ to $K_2$.
Given an $l$-ary relation $R$ in $L$ and $\baru\in V(K_1)^l$,
we have to check that $R^{K_1}\baru=R^{K_2}\phi\baru$.
The case that $\baru\in V(M)^l$ is trivial. Suppose that
$\{u_1,\ldots,u_l\}\setminus V(M)=\{w_1,\ldots,w_p\}$ is nonempty.
Note that
$\{\phi(u_1),\ldots,\phi(u_l)\}\setminus V(M)=\{\phi(w_1),\ldots,\phi(w_p)\}$
and $\{u_1,\ldots,u_l\}\cap V(M)=\{\phi(u_1),\ldots,\phi(u_l)\}\cap V(M)$.
Let $v_1,\ldots,v_p$ be pairwise distinct elements in $C$ that do not
occur in $\baru$ and hence in $\phi\baru$. Define $\psi_1$ by
$\psi_1(v_i)=w_i$ for $i\le p$ and $\psi_1(x)=x$ for all other $x$
in $V(M)$. Define $\psi_2$ similarly with the difference that
$\psi_2(v_i)=\phi(w_i)$ for $i\le p$. Obviously,
$\psi_2^{-1}\phi\baru=\psi^{-1}_1\baru$. By B2, $\psi_1$ and $\psi_2$
are partial isomorphisms from $M$ to $K_1$ and $K_2$ respectively.
Therefore
$$
R^{K_1}\baru=R^M\psi^{-1}_1\baru=R^M\psi^{-1}_2\phi\baru=R^{K_2}\phi\baru.
$$
The proof is complete.
\end{proof}

With using Definition B, the following lemma is a direct consequence
of Lemma \ref{lem:dmmpiso}.

\begin{lemma}\label{lem:dmmadd}
Let $L$ be a vocabulary with maximum relation arity $k$.
Let $M$ and $M'$ be non-isomorphic $L$-structures of orders
$n$ and $n'$ respectively and $n\le n'$. Then the bound
$$
\DD1{M,M'} < \of{1-\frac1{2k}}n+k^2-k+2
$$
may be false only if $M'=M^*\add(n'-n)v$ for some structure $M^*$
isomorphic with $M$ and $v\in V(M^*)$.
\end{lemma}

\subsubsection{An upper bound for $\dfn M$}

The following result was obtained in \cite{PVV} for graphs with the proof
easily adaptable for any structures (see Lemma 4.2 and Remark 4.9
in \cite{PVV}).

\begin{lemma}\label{lem:inpvv}
{\bf (\cite{PVV})}
Let $M$ be a structure of order $n$ with maximum relation arity $k$,
$v$ be an element of $M$ with $|[v]_M|=s\ge k$, and $M'=M\add tv$
with $t\ge1$. Then
$$
s+1\le\D{M,M'}\le\DD1{M,M'}\le s+k-1+\frac{n+1}{s+1}.
$$
\end{lemma}

Putting Lemmas \ref{lem:dmmadd} and \ref{lem:inpvv} together,
we immediately obtain an upper bound for $\D M$.
Recall that $\gs(M)=\max_{v\in V(M)}|[v]_M|$.

\begin{theorem}\label{thm:dmupper}
For a structure $M$ of order $n$ with maximum relation arity $k$,
we have
$$
\DD1M\le\max\ofcurle{\of{1-\frac1{2k}}n+k^2-k+2,\;\gs(M)+k}.
$$
\end{theorem}

\begin{proof}
Given $M$, let us summarize upper bounds we have for $\DD1{M,M'}$
for various $M'$ non-isomorphic with $M$. Denote
$$
u_{k,n}=\of{1-\frac1{2k}}n+k^2-k+2\mbox{\ \ and\ \ }
f(s)=\left\lfloor s+k-1+\frac{n+1}{s+1}\right\rfloor.
$$
If $M'=M^*\add tv$ for $M^*$ an isomorphic copy of $M$, then
\begin{equation}\label{eq:dmmgsm}
\DD1{M,M'}\le\max_{1\le s\le\gs(M)} f(s)
\end{equation}
by Lemma \ref{lem:inpvv}. Similarly,
if $M=M^*\add tv$ for $M^*$ an isomorphic copy of $M'$, then
$$
\DD1{M,M'}\le\max_{1\le s\le\gs(M')} f(s),
$$
which is within the bound \refeq{eq:dmmgsm} because in this case
$\gs(M')\le\gs(M)$. For all other $M'$ we have
$$
\DD1{M,M'} < u_{k,n}
$$
by Lemma \ref{lem:dmmadd}.

Notice now that
$$
\max_{1\le s\le\gs(M)} f(s)\le\max\{f(1),f(\gs(M))\}.
$$
Furthermore,
$$
f(\gs(M))\le\cases{
f(1) & if\ \ $\gs(M)\le(n-1)/2$,\cr
\gs(M)+k & if\ \ $\gs(M)\ge n/2$\cr}
$$
and $f(1) < u_{k,n}$.
Summing up, we conclude that
$$
\max_{M'}\DD1{M,M'}\le\max\{u_{k,n},\;\gs(M)+k\}.
$$
By Lemma \ref{lem:ddd}, the proof is complete.
\end{proof}

Note that, given $M$, the number $\gs(M)$ is efficiently computable
in the sense that computing $\gs(M)$ reduces to verification if
a transposition is an automorphism of the structure.
Thus, Theorem \ref{thm:dmupper} provides an efficiently computable
non-trivial upper bound for~$\DD1M$, whereas it seems plausible that
the exact value of $\D M$ is incomputable.

We also can restate the obtained bounds as a
dichotomy result telling us that either we have the bound
$\DD1{M,M'} < (1-\frac1{2k})n+k^2-k+2$ or else $M$ has a simple,
easily recognizable property and, moreover, for all such exceptional
$M$ we are able to easily compute $\D M$ within an additive constant.
Results of this sort are obtained in \cite{PVV} for structures
with maximum relation arity 2 and $k$-uniform hypergraphs.

\begin{theorem}
Let $M$ be a structure of order $n$ with maximum relation arity $k$. If
\begin{equation}\label{eq:dich1}
\gs(M)\le\of{1-\frac1{2k}}n+(k-1)^2+1,
\end{equation}
we have
\begin{equation}\label{eq:dich2}
\DD1{M}\le\of{1-\frac1{2k}}n+k^2-k+2.
\end{equation}
Otherwise we have
\begin{equation}\label{eq:dich3}
\gs(M)+1\le\D M\le\DD1M\le\gs(M)+k.
\end{equation}
\end{theorem}

\begin{proof}
If the condition \refeq{eq:dich1} is met, the bound \refeq{eq:dich2}
follows directly from Theorem \ref{thm:dmupper}. If \refeq{eq:dich1}
does not hold, the upper bound in \refeq{eq:dich3} again follows
from Theorem \ref{thm:dmupper}. The lower bound in \refeq{eq:dich3}
follows from Lemma \ref{lem:inpvv} as
$\D M\ge\D{M,M\add1v}\ge\gs(M)+1$, where $v\in V(M)$ is such that
$|[v]_M|=\gs(M)$ and hence $|[v]_M|>k$.
\end{proof}

\section{Identifying finite structures by Bernays-Sch\"on\-finkel formulas}%
\label{s:bs}

\begin{theorem}\label{thm:prenex}
Let $L$ be a vocabulary with maximum relation arity $k$.
If $M$ is an $L$-structure of order $n$, then
\begin{equation}
\bs M < \of{1-\frac1{2k^2+2}}n+k.
\end{equation}
If $k=1$, a stronger bound $\bs M\le n/2+1$ holds true.
\end{theorem}

The case of $k=1$ is easy and included for the sake of completeness.
The upper bound of $n/2+1$ matches, up to an additive constant of 1,
a simple lower bound of $n/2$ attainable by structures with a
single unary relation. The proof of Theorem \ref{thm:prenex} for the
case that $k\ge2$ takes the rest of this section.

\subsection{Notation}

In addition to the notation introduced in Section \ref{s:notation},
we will denote $[k]=\{1,2,\ldots,\allowbreak k\}$.
If $\bar z=(z_1,\ldots,z_l)$
and $\tau$ is a map from $[k]$ to $[l]$, then
$\bar z^\tau=(z_{\tau(1)},\ldots,z_{\tau(k)})$.

Recall that,
given a partial isomorphism $\phi\function X{X'}$ from an $L$-structure
$M$ to another $L$-structure $M'$, we have defined a relation $\phieq$
between elements in $\compl X$ and elements in $\compl{X'}$
(see Definition \ref{def:phieq}).
Definition \ref{def:phieqcl} extends this relation over classes
in $\calC(X)$ and $\calC(X')$. We will need yet another extension
of $\xeq$ over subsets of $\compl X$ and $\compl{X'}$.
Let $U\subseteq\compl X$ and $U'\subseteq\compl{X'}$.
We will write $U\phicong U'$ if $\phi$ extends to an isomorphism from
$M[X\cup U]$ to $M'[X'\cup U']$.

We define $\Bs qM$ similarly to $\bs M$ with the only additional
requirement that an identifying \BS\/ formula has at most $q$
universal quantifiers. It is clear that
$\bs M\le\Bs{q+1}M\le\Bs qM$.

\subsection{A couple of useful formulas}

If $\bar x=(x_1,\ldots,x_l)$ is a sequence of variables, let
$$
\Dist (\bar x)=\bigwedge_{1\le i<j\le l} x_i\ne x_j.
$$

Let $M$ be a finite structure over vocabulary $L$ and $\bar a$
be a sequence of $l$ pairwise distinct elements of $V(M)$.
Then it is easy to construct a first order formula
$\Iso_{M,\bar a}(x_1,\ldots,x_l)$ such that, for every $L$-structure
$M'$ and $\bar a'\in V(M')^l$, $M',\bar a'\models\Iso_{M,\bar a}(\bar x)$
iff the component-wise correspondence between $\bar a$ and $\bar a'$
is a partial isomorphism between $M$ and $M'$. Specifically, assume
that $L=(R_1,\ldots,R_m)$, where $R_i$ has arity $k_i$. Then
\begin{eqnarray*}
\Iso_{M,\bar a}(\bar x)=
\Dist(\bar x)
\And
\bigwedge_{i=1}^m
\biggl(
&\displaystyle\bigwedge_\tau&
\setdeff{R_i(\bar x^\tau)}{\tau\function{[k_i]}{[l]},\ R_i^M(\bar a^\tau)=1}\\
\And
&\displaystyle\bigwedge_\tau&
\setdeff{\neg R_i(\bar x^\tau)}{\tau\function{[k_i]}{[l]},\ R_i^M(\bar a^\tau)=0}
\biggl).
\end{eqnarray*}

\subsection{The first way of identification}

In this section we will exploit the relation $\sim$ on $V(M)$
defined in Section \ref{ss:userel} and the invariant $\gs(M)$
introduced in Definition~\ref{def:sigma}.

\begin{proposition}\label{prop:sigma}
Let $L$ be a vocabulary with maximum relation arity $k$.
For every $L$-structure $M$ of order $n$, we have
$$
\Bs kM\le n+k-\gs(M).
$$
\end{proposition}

\begin{proof}
Suppose that $\gs(M)=k+d$ with $d\ge1$ (if $\gs(M)\le k$, the proposition
is trivial). Let $A$ be a $\sim$-equivalence class of elements of $V(M)$
such that $|A|=\gs(M)$.
Denote $B=\compl A$ and fix orderings $A=\{a_1,\ldots,a_{k+d}\}$
and $B=\{b_1,\ldots,b_{n-k-d}\}$. Set $\bar a=(a_1,\ldots,a_k)$.
We suggest the following formula $\Phi_M$ to identify $M$:
$$
\Phi_M=
\exists y_1\ldots\exists y_{n-k-d}
\forall x_1\ldots\forall x_k
\Psi_M(\bar y,\bar x),
$$
where
$$
\Psi_M(\bar y,\bar x)=\Iso_{M,\bar b}(\bar y)\And
\of{
\Dist(\bar y,\bar x)\to
\Iso_{M,\bar b,\bar a}(\bar y,\bar x)
}.
$$

\begin{claim}\label{cl:s1}
$M'\models\Phi_M$ iff there is a partial isomorphism $\phi\function B{B'}$
from $M$ to $M'$ such that every injective extension of $\phi$ over
$B\cup\{a_1,\ldots,a_k\}$ is a partial isomorphism from $M$ to~$M'$.
\end{claim}

\begin{claim}\label{cl:s2}
$M\models\Phi_M$.
\end{claim}

\begin{subproof}
On the account of Claim \ref{cl:s1}, it suffices to show that the
extension $\phi$ of $\id_B$ by $\phi(a_1)=a_{i_1}$, \ldots,
$\phi(a_k)=a_{i_k}$, where $i_1,\ldots,i_k$ is an arbitrary sequence
of pairwise distinct indices in $[k+d]$, is a partial automorphism of
$M$. This follows from the fact that every permutation of $A$,
in particular, that taking each $a_j$ for $j\le k$ to $a_{i_j}$,
is decomposed into a product of transpositions $(a_p a_q)$ with
$1\le p<q\le k+d$. (Recall that the latter are automorphisms of~$M$).
\end{subproof}

\begin{claim}\label{cl:s3}
If an $L$-structure $M'$ has order $n$ and $M'\models\Phi_M$, then
$M$ and $M'$ are isomorphic.
\end{claim}

\begin{subproof}
Let $\phi$ and $B'$ be as in Claim \ref{cl:s1}. Fix an ordering
$a'_1,\ldots,a'_{k+d}$ of the set $A'=V(M')\setminus B'$.
According to Claim \ref{cl:s1}, for every sequence $i_1,\ldots,i_k$
of pairwise distinct indices in $[k+d]$, the extension of $\phi$
by $\phi(a_j)=a'_{i_j}$ for $j\le k$ is a partial isomorphism
from $M$ to $M'$. From the proof of Claim \ref{cl:s2} we know
an analog of this fact for $M$ itself:
for every sequence $i_1,\ldots,i_k$
of pairwise distinct indices in $[k+d]$, the extension $\psi$ of $\id_B$
by $\psi(a_j)=a_{i_j}$ for $j\le k$ is a partial automorphism
of $M$. It follows, in particular, that for every sequence
$1\le i_1<\ldots<i_k\le k+d$, the extension of $\phi$ by
$\phi(a_{i_j})=a'_{i_j}$ for $j\le k$ is a partial isomorphism
from $M$ to $M'$. Extend $\phi$ over the whole $V(M)$ by
$\phi(a_i)=a'_i$ for all $i\le k+d$. We conclude that the restriction
of $\phi$ on every $k$-element subset of $V(M)$ is a partial isomorphism
from $M$ to $M'$. Since every relation of $M$ has arity at most $k$,
$\phi$ is an isomorphism from $M$ to~$M'$.
\end{subproof}
\end{proof}

\subsection{The second way of identification}

\begin{definition}\label{def:base2}
A set $B\subseteq V(M)$ is called a {\em base\/} of a structure $M$
if the relations $\beq$ and $\sim$ coincide on $\compl B$.
The {\em fineness\/} of a base $B$ is defined by
$f(B)=\max\setdef{|C|}{C\in\calC(B)}$. Furthermore, let
$\rho(B)=|B|+\max\{f(B)+1,k\}$.

We define $\rho(M)$ to be the minimum $\rho(B)$ over all bases $B$ of~$M$.
\end{definition}

\begin{proposition}\label{prop:rho}
$\bs M\le\rho(M)$.
\end{proposition}

\begin{proof}
Given a base $B$ of $M$, we construct a \BS\/ formula $\Phi_M$ with
$\rho(B)$ quantifiers that identifies $M$. Let $p=|B|$ and
$q=\max\{f(B)+1,k\}$. Assume that $p+q<n$ for otherwise we are done.
Denote $A=\compl B$ and fix orderings
$A=\{a_1,\ldots,a_{n-p}\}$ and $B=\{b_1,\ldots,b_p\}$. We set
$$
\Phi_M=
\exists y_1\ldots\exists y_p
\forall x_1\ldots\forall x_q
\Psi_M(\bar y,\bar x),
$$
where
$$
\Psi_M(\bar y,\bar x)=\Iso_{M,\bar b}(\bar y)\And
\Biggl(
\Dist(\bar y,\bar x)\to
\bigvee_{
\begin{array}{c}
\scriptstyle\tau\function{[q]}{[n-p]}\\
\scriptstyle\tau\mbox{\scriptsize\ is\ injective}
\end{array}
}
\Iso_{M,\bar b,\bar a^\tau}(\bar y,\bar x)
\Biggr).
$$

\begin{claim}\label{cl:r1}
Let $M'$ be another $L$-structure, $\bar b'=(b'_1,\ldots,b'_p)$
be a sequence of elements of $V(M')$, and
$A'=V(M')\setminus\{b'_1,\ldots,b'_p\}$. Then
$M',\bar b'\models\forall x_1\ldots\forall x_q\Psi_M(\bar y,\bar x)$
holds iff
\begin{itemize}
\item
the component-wise correspondence $\phi$ between $\bar b$
and $\bar b'$ is a partial isomorphism from $M$ to $M'$ and
\item
for every $U'\subseteq A'$ with at most $q$ elements there is a
$U\subseteq A$ such that $U\phicong U'$.
\end{itemize}
\end{claim}

The proof is fairly obvious. The claim immediately implies that
$M\models\Phi_M$.

\begin{claim}\label{cl:r2}
If $M'\models\Phi_M$ and $M'$ has order $n$, then $M$ and $M'$ are
isomorphic.
\end{claim}

\begin{subproof}
Let $\bar b'=(b'_1,\ldots,b'_p)$ be such that
$$
M',\bar b'\models\forall x_1\ldots\forall x_q\Psi_M(\bar y,\bar x).
$$
Set $B'=\{b'_1,\ldots,b'_p\}$. By the definition of $\Psi_M$,
there is a partial isomorphism $\phi\function B{B'}$ from $M$ to $M'$.
By Claim \ref{cl:r1}, for every $a'\in A'$ there is $a\in A$ such that
$a\phieq a'$. Hence for every $C'\in\calC(B')$ there is $C\in\calC(B)$
such that $C\phieq C'$. Moreover, for every $C'\in\calC(B')$
and the respective $C\in\calC(B)$ it holds $|C|\ge|C'|$
(if $|C'|>|C|$, then for any $(|C|+1)$-element set $U'\subseteq C'$
the second condition in Claim \ref{cl:r1} fails). Since $|A|=|A'|$
or, in other terms, $\sum_{C\in\calC(B)}|C|=\sum_{C'\in\calC(B')}|C'|$,
for every $C'$ it actually holds the equality $|C|=|C'|$.
Thus, we have a one-to-one correspondence between $\calC(B)$
and $\calC(B')$ such that, if $C\in\calC(B)$ and $C'\in\calC(B')$
correspond to one another, then $C\phieq C'$ and $|C|=|C'|$.

We are now prepared to exhibit an isomorphism from $M'$ to $M$.
Fix an arbitrary extension $\psi$ of $\phi^{-1}$ to a one-to-one
map from $V(M')$ to $V(M)$ taking each $C'$ to the respective $C$.
We will show that $\psi$ is an isomorphism. Let $R'$ be an $l$-ary
relation of $M'$ and $R$ be the respective relation of $M$.
Given an arbitrary $l$-tuple $\bar u'\in V(M')^l$, we have to prove
that
\begin{equation}\label{eq:R=R}
R\psi\bar u'=R'\bar u'.
\end{equation}

Denote $U'=\{u'_1,\ldots,u'_l\}$. Let $\psi_{U'}$ be the extension of
$\phi^{-1}$ to a partial isomorphism from $M'$ to $M$ with
$U'\subseteq\dom\psi_{U'}$ whose existence is guaranteed by
Claim \ref{cl:r1}. We have
$$
R\psi_{U'}\bar u'=R'\bar u'.
$$
To prove \refeq{eq:R=R}, it suffices to prove that
\begin{equation}\label{eq:R=R:U'}
R\psi_{U'}\bar u'=R\psi\bar u'.
\end{equation}

We proceed similarly to the proof of Claim \ref{cl:upsilon}
in Section \ref{ss:strategy}.
By Item 3 of Lemma \ref{lem:phieq}, the partial map $\psi_{U'}$
takes an element in a class $C'$ to an element in the respective class $C$.
Suppose that $\psi_{U'}$ is extended
over the whole $V(M')$ with the latter condition obeyed.
Since both $\psi_{U'}$ and $\psi$
extend $\phi^{-1}$, the product $\psi_{U'}\psi^{-1}$ moves
only elements in $A$.
Since both $\psi$ and $\psi_{U'}$
take an element in a class $C'$ to an element in the respective class $C$,
the map $\psi_{U'}\psi^{-1}$ preserves the partition
$\calC(B)$ of $A$.
It follows that $\psi_{U'}\psi^{-1}$ is decomposed into the product
of permutations $\pi_C$ over $C\in\calC(B)$, where each $\pi_C$ acts
on the respective $C$. Since every $\pi_C$ is decomposable into
a product of transpositions, we have
$\psi_{U'}\psi^{-1}=\tau_1\tau_2\ldots\tau_t$
with $\tau_i$ being a transposition of two elements both in some $C$.
It is easy to see that
$\psi_{U'}\bar u'=(\ldots((\psi\bar u')^{\tau_t})\ldots)^{\tau_1}$.
By Lemma \ref{lem:xeqsim}, each application of $\tau_i$ does not change
the initial value of $R\psi\bar a'$. Therewith \refeq{eq:R=R:U'}
is proved.
\end{subproof}
\end{proof}

\begin{remark}\label{rem:rhodef}
One can show that $\rho(M)$ provides us with an upper bound not only
for $\bs M$ but also for $\DD1M$.
\end{remark}

\subsection{The third way of identification}

Yet another way of identification that we here suggest is actually
not new being a specification of Proposition \ref{prop:rho} in the
preceding section.

\begin{definition}
If $M$ is a finite structure, let
$$
\gd(M)=\max\setdef{|A|}{A\subseteq V(M)\mbox{\ such\ that\ }
a_1\notaeq a_2\mbox{\ for\ every\ }a_1,a_2\in A}.
$$
\end{definition}
It is not hard to see that, in other terms,
$\gd(M)=\max_{X\subseteq V(M)}|\calC(X)|$.

\begin{proposition}\label{prop:delta}
Let $L$ be a vocabulary with maximum relation arity $k\ge2$.
For every $L$-structure $M$ of order $n$, we have
$$
\Bs kM\le n+k-\gd(M).
$$
\end{proposition}

\begin{proof}
As easily seen, if $A\subseteq V(M)$ is such that
$a_1\notaeq a_2$ for every $a_1,a_2\in A$, then $\compl A=V(M)\setminus A$
is a base of $M$ with fineness $f(\compl A)=1$. Since $k\ge2$, we
have $\max\{f(\compl A)+1,k\}=k$ and therefore
$\rho(M)\le n+k-\gd(M)$. Thus, the proposition directly follows from
Proposition~\ref{prop:rho}. We only have to note that the
identifying formula constructed in the proof of Proposition~\ref{prop:rho}
has $\max\{f(\compl A)+1,k\}=k$ universal quantifiers.
\end{proof}

\subsection{Putting it together}

We now complete the proof of Theorem \ref{thm:prenex}.
Assume that $k\ge2$.
We will employ all three possibilities of identifying $M$ given by
Propositions \ref{prop:delta}, \ref{prop:sigma}, and \ref{prop:rho}.
Using the last possibility, we will use the set $X_{k+1}$
defined by Definition \ref{def:base} that is a base of $M$ according to
Lemma \ref{lem:xeqsim}.

By the bound \refeq{eq:a0} of Lemma \ref{lem:unpebbl} and the fact
that $|\calC(X)|\le\gd(M)$ for every $X\subseteq V(M)$, we have
\begin{equation}\label{eq:X=n_le}
|X_{k+1}|=n-|Z|\le2k^2\gd(M)-(k-1).
\end{equation}
We now consider two cases.

\case 1{$Z=\emptyset$.}
By \refeq{eq:X=n_le} we have $\gd(M)\ge\frac{n+k-1}{2k^2}$.
By Proposition \ref{prop:delta}, this implies that
$$
\bs M < \of{1-\frac1{2k^2}}n+k.
$$

\case 2{$Z\ne\emptyset$.}
In this case for the fineness of the base $X_{k+1}$ we have
$f(X_{k+1})\ge k+2$. Using \refeq{eq:X=n_le}, we obtain
$$
\rho(X_{k+1})\le2k^2\gd(M)-(k-1)+\max\setdef{|C|}{C\in\calC(X_{k+1})}+1\le
2k^2\gd(M)+\gs(M)+2-k.
$$
Let $\lambda(M)=\max\{\gd(M),\gs(M)\}$. By Propositions
\ref{prop:delta}, \ref{prop:sigma}, and \ref{prop:rho}, we have
\begin{eqnarray*}
\bs M&\le&\min\{n+k-\gd(M),n+k-\gs(M),2k^2\gd(M)+\gs(M)+2-k\}\\&\le&
\min\{n+k-\lambda(M),(2k^2+1)\lambda(M)+2-k\}\\
&\le&\max_{1\le\lambda\le n}\min\{n+k-\lambda,(2k^2+1)\lambda+2-k\}\\&\le&
\of{1-\frac1{2k^2+2}}n+k-\frac{k-1}{k^2+1}<\of{1-\frac1{2k^2+2}}n+k.
\end{eqnarray*}
Since the latter bound holds in both the cases, the proof
of Theorem \ref{thm:prenex} is for $k\ge2$ complete.

In the case of $k=1$ we use Propositions \ref{prop:sigma} and \ref{prop:rho}.
We use the fact that, for a structure $M$ with all relations unary,
the empty set is a base and $\rho(\emptyset)=\gs(M)+1$. We therefore have
$\bs M\le\min\{n+1-\gs(M),\gs(M)+1\}\le n/2+1$.

\section{Identifying finite structures by Bernays-Sch\"on\-finkel formulas
with bounded number of universal quantifiers}\label{s:unibound}

Recall that
$\Bs qM$ denotes the minimum total number of quantifiers in a
\BS\/ formula identifying $M$ with at most $q$ universal quantifiers.
We now address the asymptotics of the maximum value of $\Bs qM$
over structures of order $n$ under the condition that $q$
is bounded by a constant.
We first observe that less than $k$ universal quantifiers
are rather useless for identification of a structure with
maximum relation arity~$k$.

\begin{proposition}
If $M$ is a structure of order $n$ with maximum relation arity $k$
and $n\ge k$, then $\Bs{k-1}M=n$.
\end{proposition}

\begin{proof}
We have to show that no formula
$\Phi=\exists y_1\ldots\exists y_p\forall x_1\ldots\forall x_q
\Psi(\bar y,\bar x)$
with $\Psi$ quantifier-free, $q\le k-1$, and $p+q\le n-1$ can identify
$M$. Suppose that
\begin{equation}\label{eq:mba_mod_psi}
M,\bar b,\bar a\models\Psi(\bar y,\bar x)
\end{equation}
for some $\bar b\in V(M)^p$ and all $\bar a\in V(M)^q$.
Let $A=V(M)\setminus\{b_1,\ldots,b_p\}$. Since $q+1\le k$,
$q+1\le n-p\le|A|$, and $n\ge k$, there is a $k$-element $U\subseteq V(M)$
such that $|U\cap A|\ge q+1$. Let $u_1,\ldots,u_k$ be an arbitrary
ordering of $U$. Let $R$ be a $k$-ary relation of $M$.
Define a relation $R'$ so that $R'\bar u\ne R\bar u$ and $R'$
coincides with $R$ elsewhere. Let $M'$ be the modification of $M$
with $R'$ instead of $R$. Clearly, $M'$ and $M$ are non-isomorphic.
It is easy to see that $M',\bar b,\bar a\models\Psi(\bar y,\bar x)$
for the same $\bar b$ as in \refeq{eq:mba_mod_psi} and all
$\bar a\in V(M')^q$. Therefore $M'\models\Phi$ and $\Phi$ fails
to identify~$M$.
\end{proof}

If at least $k$ universal quantifiers are available, some saving
on the number of quantifiers is possible: It turns out that
$\Bs kM < n-\sqrt n+k^2+k$ and this bound cannot be improved much
if we keep the number of universal quantifiers constant.

\begin{theorem}\label{thm:bsunib}
Let $\Bs q{n,k}$ denote the maximum $\Bs qM$ over structures $M$
of order $n$ and maximum relation arity $k$. Then
$$
\Bs k{n,k} < n-\sqrt n+k^2+k.
$$
On the other hand, if $n$ is a square, then
$$
\Bs q{n,k}\ge n-(q-1)\sqrt n+q
$$
for every $q\ge 2$ and $k\ge 2$.
\end{theorem}

The upper bound of Theorem \ref{thm:bsunib} is provable by the
techniques from Section \ref{s:bs}.
Let $M$ be a structure of order $n$ with maximum relation arity $k$.
By Propositions \ref{prop:delta} and \ref{prop:sigma},
$$
\Bs kM\le n+k-\max\{\gd(M),\gs(M)\}.
$$
It remains to prove the following bound.

\begin{lemma}\label{lem:maxds}
$\max\{\gd(M),\gs(M)\}>\sqrt n-k^2$.
\end{lemma}

\begin{proof}
By the bound \refeq{eq:a0} of Lemma \ref{lem:unpebbl},
$$
n+k-1\le 2k\sum^{k-1}_{i=1}|\calC^{k+1}(X_i)|+
(k+1)|\calC^{k+1}(X_k)|+(k-1)|\calC(X_k)|+|Z|.
$$
We bound each term $|\calC(X)|$ from above by $\gd(M)$.
Furthermore, we bound $|Z|$ from above by the number of
$\xkkeq$-equivalence classes inside $Z$ multiplied by
the maximum number of elements in such a class. By Lemma \ref{lem:xeqsim}
it follows that $|Z|\le \gd(M)\gs(M)$. We therefore conclude that
$$
n+k-1\le\gd(M)(2k^2+\gs(M)).
$$
This implies
$$
\max\{\gd(M),\gs(M)\}\ge
\min_{1\le\gs\le n}\max\ofcurle{\gs,\frac{n+k-1}{2k^2+\gs}}>
\sqrt n-k^2,
$$
as required.
\end{proof}

\begin{remark}\label{rem:maxdsopt}
The bound of Lemma \ref{lem:maxds} is essentially optimal because,
for any graph $G$ of order $m^2$ whose vertex set is partitioned
into $m$ $\sim$-equivalence classes of $m$ element each, it holds
$\gs(G)=m$ and $\gd(G)\le m$. Such $G$ can be constructed from any
graph $H$ of order $m$ whose automorphism group contains no transposition
by replacing each vertex $v\in V(H)$ with $m$ pairwise (non-)adjacent
vertices $\sim$-related to $v$ in~$H$.
\end{remark}

We now prove the lower bound of Theorem \ref{thm:bsunib}.
It suffices to do it for graphs. The example of $G$ with large $\Bs qG$
will be the same as in Remark \ref{rem:maxdsopt}. This example
can be lifted to a higher arity $k$ by adding $k-2$ dummy coordinates
to the adjacency relation with no affect to its truth value.

\begin{proposition}
Let $G_m$ be graph of order $m^2$ whose vertex set is partitioned
into $m$ $\sim$-equivalence classes of $m$ element each. Let $q\ge 2$.
Then $\Bs q{G_m}\ge m^2-(q-1)m+q$.
\end{proposition}

\begin{proof}
It is enough to show that, if $G_m$ is identified by a \BS\/ formula $\Phi$
with $q$ universal quantifiers, then $\Phi$ contains at least
$m^2-(q-1)m$ existential quantifiers. If $q\ge m+1$, this is trivial.
Assume that $q\le m$.

Suppose on the contrary that $G_m$ is
identified by a \BS\/ formula $\Phi=\exists y_1\ldots\exists y_p
\forall x_1\ldots\forall x_q \Psi(\bar y,\bar x)$
with $p<m^2-(q-1)m$. Let $\bar b\in V(G_m)^p$ be such that
$G_m,\bar b\models\forall x_1\ldots\forall x_q \Psi(\bar y,\bar x)$.
Equivalently,
\begin{equation}\label{eq:gmba_mod_psi}
G_m,\bar b,\bar a\models \Psi(\bar y,\bar x)
\mbox{\ for\ all\ }\bar a\in V(G_m)^q.
\end{equation}
Let $A=V(G_m)\setminus\{b_1,\ldots,b_p\}$. We have $|A|\ge(q-1)m+1$.
The condition imposed on $G_m$ implies that there are two
$\sim$-equivalence classes, $C_1$ and $C_2$, such that
$|A\cap C_1|\ge q$ and $|A\cap C_2|\ge1$. Let us modify $G_m$
by removing one vertex from $A\cap C_2$ and adding a new vertex $v'$
to $C_1$ so that $v'\sim v$ for all $v\in C_1$. The modified graph, $G'$,
is clearly non-isomorphic to $G_m$. We show that, nevertheless,
$G'\models\Phi$.

It suffices to show that $G',\bar b,\bar a'\models\Psi(\bar y,\bar x)$
for every $\bar a'\in V(G')^q$. In view of \refeq{eq:gmba_mod_psi},
we are done if for every $\bar a'\in V(G')^q$ we are able to find
an $\bar a\in V(G_m)^q$ such that the component-wise correspondence
between $\bar b,\bar a$ and $\bar b,\bar a'$ is a partial isomorphism
between $G_m$ and $G'$. If $\bar a'$ does not contain any occurrence of $v'$,
we obviously can take $\bar a=\bar a'$. If $\bar a'$ contains an
occurrence of $v'$, let $v$ be a vertex in $A\cap C_1$ that does not
occur in $\bar a'$ and let $\bar a$ be the result of substituting $v$
in place of $v'$ everywhere in $\bar a'$. It is not hard to see that
the obtained $\bar a$ is as required.
\end{proof}

\section{The case of graphs}\label{s:graphs}

For a binary structure $M$, Theorem \ref{thm:main} implies
$\I M<0.75n+4$ and Theorem \ref{thm:prenex} implies $\bs M<0.9n+2$.
In the case of graphs, both these bounds can be improved.
In \cite{PVV} we obtain an almost optimal bound $\I G\le (n+3)/2$
(there are simple examples of graphs with $\I G\ge (n+1)/2$).
Combining the approach from \cite{PVV} and the techniques from
Section \ref{s:bs}, we are able to prove a better bound for $\bs G$
as well. We are also interested
in knowing the smallest $n$ starting from which for $G$ of order $n$
we have the bound at least $\bs G\le n-1$, an improvement on the trivial
bound of $n$.

\begin{theorem}\label{thm:graphs}
Let $G$ be a graph of order $n$.
\begin{enumerate}
\item
We have $\bs G\le 3n/4+3/2$.
\item
If $n\ge 5$, we have $\Bs 2G\le n-1$ with the only exception
of the graph $H$ on 5 vertices with 2 adjacent edges for which,
nevertheless, we have $\Bs 3H\le4$.
\end{enumerate}
\end{theorem}

\begin{proof}
Given a graph $G$,
let $X=E(\emptyset)$, where the transformation $E$ is introduced in
Section \ref{ss:transform}. We state two properties of the $X$
established in \cite{PVV}:
\begin{description}
\item[{\em Property 1.}]
$|\calC(X)|\ge|X|+1$.
\item[{\em Property 2.}]
Let $Y=Y(X)$ and $Z=V(G)\setminus(X\cup Y)$.
Every class in $\calC(X\cup Y)$ consists of pairwise $\sim$-equivalent
vertices.
\end{description}
Note that $|\calC(X)|\le\gd(G)$. By Property 1 we conclude that
\begin{equation}\label{eq:gr3}
|X|+|Y|\le|X|+|\calC(X)|\le2|\calC(X)|-1\le2\gd(G)-1.
\end{equation}
Property 2 means that $X\cup Y$ is a base of $G$.

We now consider two cases.

\case 1{$Z=\emptyset$.}
In this case $n=|X|+|Y|$. By \refeq{eq:gr3}, $\gd(G)\ge(n+1)/2$.
Using Proposition \ref{prop:delta}, we obtain $\Bs 2G\le n/2+3/2$.

\case 2{$Z\ne\emptyset$.}
In this case $\rho(X\cup Y)\le |X|+|Y|+\gs(G)+1$. By \refeq{eq:gr3}
we have $\rho(X\cup Y)\le 2\gd(G)+\gs(G)$.
Denote $\lambda(G)=\max\{\gd(G),\gs(G)\}$. By Propositions
\ref{prop:delta}, \ref{prop:sigma}, and \ref{prop:rho}, we have
\begin{eqnarray*}
\bs G&\le&\min\{n+2-\gd(G),n+2-\gs(G),2\gd(G)+\gs(G)\}\\&\le&
\min\{n+2-\lambda(G),3\lambda(G)\}\\&\le&
\max_{1\le\lambda\le n}\min\{n+2-\lambda,3\lambda\}=3n/4+3/2.
\end{eqnarray*}
Since this bound holds true in both the cases, Item 1 of the theorem
is proved.

To prove Item 2, we estimate $\max\{\gd(G),\gs(G)\}$.
Since $n=|X|+|Y|+|Z|\le 2\gd(G)-1+\gd(G)\gs(G)$, we have
\begin{equation}\label{eq:seven}
n+1\le\gd(G)(2+\gs(G)).
\end{equation}
It follows that
$$
\max\{\gd(G),\gs(G)\}\ge\min_{1\le c\le n}\max\ofcurle{c,\frac{n+1}{2+c}}=
\sqrt{n+2}-1
$$
and hence
\begin{equation}\label{eq:max3}
\max\{\gd(G),\gs(G)\}\ge 3
\end{equation}
whenever $n>7$.

\begin{claim}\label{cl:six}
The bound \refeq{eq:max3} holds for all $G$ of order~6.
\end{claim}

\noindent
This claim is proved by the direct brute force analysis.
Making it, it suffices to consider graphs on 6 vertices with
at most 7 edges. The reason is that $\gd(\compl G)=\gd(G)$
and $\gs(\compl G)=\gs(G)$, where $\compl G$ denotes the complement
of a graph $G$, i.e., the graph with the same vertex set and
exactly those edges absent in~$G$.

\begin{claim}
The bound \refeq{eq:max3} holds for all $G$ of order~7.
\end{claim}

\begin{subproof}
Given three vertices $x$, $y$, and $z$, we say that $x$ separates
$y$ and $z$ if $x$ is adjacent to exactly one of $y$ and $z$.
Note that $y\sim z$ iff no $x$ separates these vertices.

Let a graph $G$ have 7 vertices. If $\gs(G)=1$, then $\gd(G)\ge 3$ by
\refeq{eq:seven}. Our task is therefore to deduce $\gd(G)\ge3$
from the assumption that $\gs(G)=2$.

Let $u$ and $v$ be $\sim$-equivalent vertices of $G$. Suppose
that they are adjacent. We do not lose the generality because
it suffices to prove the claim for one of $G$ or $\compl G$.
Let us remove $u$, that is, consider the graph $G-u=G[V(G)\setminus\{u\}]$.
As it is easy to see, $\gd(G)\ge \gd(G-u)$. Thus, if
$\gd(G-u)\ge 3$, we are done. Otherwise, by Claim \ref{cl:six},
we have $\gs(G-u)\ge 3$.

Let $a$, $b$, and $c$ be pairwise $\sim$-equivalent vertices in $G-u$.
Assume for a while that $v\notin\{a,b,c\}$. Since $u\sim v$ in $G$,
the vertex $u$ does not separate $a$, $b$, and $c$ because $v$ does not.
Thus, these three vertices are pairwise $\sim$-equivalent in $G$,
contradicting our assumption that $\gs(G)=2$. We conclude that
$v\in\{a,b,c\}$.

Without loss of generality, assume that $v=c$.
Since $v$ cannot be $\sim$-equivalent with $a$ or $b$ in $G$
and only $u$ can separate $v$ from $a$ and from $b$,
the vertex $u$ must be non-adjacent to $a$ and to $b$.
The same is true for the $\sim$-equivalent vertex $v$.
Thus, $\{u,a\}$, $\{u,b\}$, $\{v,a\}$, and $\{v,b\}$ all are non-edges.
Since $\{a,v\}$ is a non-edge, $\{a,b\}$ is a non-edge too
because $v=c\sim b$ in $G-u$. Note that $a\sim b$ in $G$ because
this is so in $G-u$ and $u$ does not separate these two vertices.

Apply now the same trick with removal of $a$ instead of $u$. If we are not
done, then $G-a$ has a $\sim$-equivalence class $\{b,s,t\}$.
Our argument is completely symmetric with the difference that
$a$ and $b$ are now non-adjacent vertices. We hence should switch over
all (non)adjacencies and conclude that
$\{a,s\}$, $\{a,t\}$, $\{b,s\}$, $\{b,t\}$, and $\{s,t\}$ all
are edges of $G$. Furhtermore, $s\sim t$ in~$G$.

Notice that $s$, $t$, $u$, and $v$ are pairwise distinct.
Indeed, since $u$ separates $b$ and $v$, we have $v\notin\{s,t\}$.
Similarly, $u\notin\{s,t\}$.

We apply the same trick once again, now with removal of $s$.
As above, unless we are done, $G-s$ has a $\sim$-equivalence
class $A$ containing $t$ and
at least 2 more vertices. As above, the vertices $a$ and $b$ cannot
belong to $A$. It follows that $A$ contains at least one of $u$ and $v$.
Actually $A$ must contain both $u$ and $v$ because these vertices
are $\sim$-equivalent.
As $A$ is not a $\sim$-class in $G$, the vertex $s$ must separate
$t$ from $u$ and from $v$. Therefore $\{s,u\}$ and $\{s,v\}$
are non-edges. As $t\sim s$, $\{t,u\}$ and $\{t,v\}$ are non-edges
too. But now the triple $a,s,u$ shows
that $\delta(G)\ge 3$. Indeed, $v$ separates $u$ from $a$ and from $s$
while $b$ separates $a$ and $s$.
\end{subproof}

Thus, if $G$ has order at least 6, we have the bound \refeq{eq:max3} and
the theorem follows from Propositions \ref{prop:delta} and \ref{prop:sigma}.
For graphs of order 5 the estimate \refeq{eq:max3} holds with the only
exception for the specified graph $H$. This graph is identified by formula
$$
\exists y_1\forall x_1\forall x_2\forall x_3
\of{\Dist(y_1,x_1,x_2,x_3)\to
\neg E(x_1,x_2)
\And
\bigvee_{i=1}^3 E(y_1,x_i)
\And
\bigvee_{i=1}^3 \neg E(y_1,x_i)
},
$$
where $E$ is the adjacency relation.
\end{proof}

\begin{remark}
Item 2 of Theorem \ref{thm:graphs} does not hold true for graphs
of order $n=4$: It is not hard to prove that $\bs F=4$ for the graph
$F$ on 4 vertices with 1 edge.
\end{remark}

\begin{remark}
In \cite{KPSV} we address the first order definability of a random
graph $G$ on $n$ vertices. It is proved that, with probability $1-o(1)$,
$$
\log_2n-2\log_2\log_2n\le \I G
\le\log_2n-\log_2\log_2n+\log_2\log_2\log_2n+O(1).
$$
One of the ingredients of the proof is that, with high probability,
$\gd(G)\ge n-(2+o(1))\log_2n$. Since $\I G\le\Bs2G\le n+2-\gd(G)$,
we conclude that, with high probability,
$$
\log_2n-2\log_2\log_2n\le \Bs2G\le (2+o(1))\log_2n.
$$
\end{remark}

\section{Open problems}\label{s:open}
\mbox{}

\que
Let $\I{n,k}$ (resp.\ $\II l{n,k}$; $\bs{n,k}$) be the maximum
$\I{M}$ (resp.\ $\II l{M}$; $\bs{M}$) over structures of order $n$
with maximum relation arity $k$. We now know that
\begin{equation}\label{eq:now}
\begin{array}{rcl}
\frac n2\le\I{n,k}\le&\ant\II1{n,k}&\ant\le\bs{n,k}
<(1-\frac1{2k^2+2})n+k\\[2mm]
\mbox{and\ }&\ant\II1{n,k}&\ant<(1-\frac1{2k})n+k^2-k+2.
\end{array}
\end{equation}
Note that $\I{n,k}\le\I{n,k+1}$ and that the lower bound of $n/2$ is
actually for $\I{n,1}$.
Make the gap between the lower and upper bounds
in \refeq{eq:now} closer.

The case of $k=2$ is essentially solved in \cite{PVV},
where the bounds
$$
\frac{n+1}2\le\I{n,2}\le\II1{n,2}\le\frac{n+3}2
$$
are proved. If $k=3$, we are able to improve on \refeq{eq:now}
by showing that $\II1{n,3}\le\frac23\,n+O(1)$ (in \cite{PVV} this bound
was obtained for 3-uniform hypergraphs).

\que
Can one prove a non-trivial upper bound for $\II0{n,k}$?
The weakest question is if $\II0{n,k}\le n-1$. It is easy to show
that $\II0{n,1}\le(n+1)/2$. In \cite{PVV} we prove that
$\II0G\le(n+5)/2$ for graphs of order~$n$.

\que
What happens if we restrict the number of existential rather than
universal quantifiers in an identifying \BS\/ formula?

\end{document}